\documentclass{article}

\usepackage[english]{babel}

\usepackage[a4paper,top=2.5cm,bottom=2.5cm,left=3cm,right=3cm,marginparwidth=1.75cm]{geometry}

\usepackage{amsmath}
\usepackage{mathtools}
\usepackage{graphicx}
\usepackage{xcolor}
\usepackage[hidelinks]{hyperref}
\usepackage{amssymb}
\usepackage{amsthm}
\usepackage{subcaption}
\usepackage{tikz}
\usepackage{colortbl}
\usepackage{makecell}

\usepackage{booktabs,multirow,array}

\usepackage{pdflscape}
\usepackage{afterpage}

\usepackage{enumitem}

\usepackage{multirow}

\usepackage[normalem]{ulem}

\graphicspath{{./fig/}}

\DeclareMathOperator*{\conv}{conv}
\DeclareMathOperator*{\ri}{ri}

\usepackage[vlined,boxed,linesnumbered]{algorithm2e}
\SetAlgoVlined
\SetKwInput{Input}{input}
\SetKwInput{Inputs}{inputs}
\SetKw{Procedure}{procedure}

\newtheorem{theorem}{Theorem}
\newtheorem{lemma}{Lemma}
\newtheorem{corollary}{Corollary}
\newtheorem{proposition}{Proposition}
\newtheorem{example}{Example}
\newtheorem{remark}{Remark}

\usepackage{stmaryrd} 
\newcommand{\iverson}[1]{\llbracket#1\rrbracket}

\newcommand*\circled[1]{\tikz[baseline=(char.base)]{%
            \node[shape=circle,draw,inner sep=2pt] (char) {#1};}}

\usepackage{enumitem}

\newcommand*{\V}{\mathcal{V}}
\newcommand*{\E}{\mathcal{E}}
\newcommand*{\F}{\mathcal{F}}
\newcommand*{\N}{\mathcal{N}}
\renewcommand*{\P}{\mathcal{P}}
\renewcommand*{\k}{k}
\renewcommand*{\l}{\ell}
\renewcommand*{\L}{\mathcal{L}}

\makeatletter
    \@ifdefinable{\bcaalg}{\def\bcaalg/{\texttt{BCA}}}
    \@ifdefinable{\hungalg}{\def\hungalg/{\texttt{Hung}}}
    \@ifdefinable{\hungrialg}{\def\hungrialg/{\texttt{Hung+RI}}}
    \@ifdefinable{\gurobialg}{\def\gurobialg/{\texttt{Gurobi}}}
\makeatother

\long\def\val#1{{\normalfont{\textsf{#1}}}}

\long\def\instancegroup#1{\emph{#1}}

\usepackage{calc}  
\newlength\mylen
\newcolumntype{P}{>{\centering\arraybackslash}p{\mylen}}

\title{Relative-Interior Solution for the (Incomplete) Linear Assignment Problem with Applications to the Quadratic Assignment Problem}
\author{Tom\'{a}\v{s} Dlask$^1$ and Bogdan Savchynskyy$^2$}
\date{$^1$ Faculty of Electrical Engineering, Czech Technical University in Prague\\ [1ex]
      $^2$ Computer Vision and Learning Lab, IWR, University of Heidelberg \\ [2ex]
      {\normalsize \texttt{dlaskto2@fel.cvut.cz, bogdan.savchynskyy@iwr.uni-heidelberg.de}}}

\usepackage[absolute,overlay]{textpos}

\begin{document}

\begin{textblock*}{21cm}(0cm,0.3cm)
\centering
The version of record of this article is published in Annals of Mathematics and Artificial Intelligence and is available online at \url{https://doi.org/10.1007/s10472-025-09974-w}.
\end{textblock*}

\maketitle

\begin{abstract}
We study the set of optimal solutions of the dual linear programming formulation of the linear assignment problem (LAP) to propose a method for computing a solution from the relative interior of this set. Assuming that an arbitrary dual-optimal solution and an optimal assignment are available (for which many efficient algorithms already exist), our method computes a relative-interior solution in linear time. Since the LAP occurs as a subproblem in the linear programming (LP) relaxation of the quadratic assignment problem (QAP), we employ our method as a new component in the family of dual-ascent algorithms that provide bounds on the optimal value of the QAP. To make our results applicable to the incomplete QAP, which is of interest in practical use-cases, we also provide a linear-time reduction from the incomplete LAP to the complete LAP along with a mapping that preserves optimality and membership in the relative interior. Our experiments on publicly available benchmarks indicate that our approach with relative-interior solution can frequently provide bounds near the optimum of the LP relaxation and its runtime is much lower when compared to a commercial LP solver.
\end{abstract}

\section{Introduction}

The NP-hard \emph{quadratic assignment problem} (QAP) is a well-studied problem of combinatorial optimization with many real-world applications, such as facility location, scheduling, data analysis, ergonomic design, or problems originating in computer vision~\cite{burkard2009assignment,cela2013quadratic,haller2022comparative}. Informally, the QAP seeks to find a bijection between two finite sets of equal size that minimizes the objective which is a sum of unary and binary functions that depend on the values of the bijection. This problem becomes the polynomially solvable \emph{linear assignment problem} (LAP) if the objective contains only unary functions (i.e., the binary functions are identically zero). There exist many algorithms for solving the QAP exactly and heuristics for providing good solutions that are not guaranteed to be optimal. For a comprehensive introduction to the LAP and the QAP, we refer to~\cite{burkard2009assignment,cela2013quadratic}. 

Importantly, our focus significantly differs from the focus of classical operations research works: the QAP instances we consider are
\begin{itemize}
\item incomplete (see definition below and Section~\ref{se:iQAP}) -- this problem class has not been considered in operations research literature so far and is mainly utilized in applied research areas like computer vision and machine learning,
\item large-scale, up to $500\times 1000$ in size,
\item extremely sparse, some with only $0.1$\% of non-zero costs, and
\item with general Lawler type of costs~\cite{lawler1963quadratic}.
\end{itemize}
For comparison, the classical operations research setting represented by the QAPLIB benchmark~\cite{burkard1997qaplib} contains complete QAP instances with sizes up to $100\times 100$, frequently much denser, and with Koopmans-Beckmann type of costs~\cite{koopmans1957assignment}. The shear size of the problems makes many powerful techniques such as SDP relaxation, interior and simplex methods very slow. Thus, we restrict our further attention only to the works that have shown their efficiency in the setting described above~\cite{haller2022comparative} and are most relevant to our approach.

\paragraph{Related Work}
We follow up on two methods for obtaining bounds on the QAP, namely~\cite{zhang2016pairwise} and~\cite{hutschenreiter2021fusion}. Both of these methods are based on \emph{block-coordinate ascent}\footnote{Formally, the QAP is defined in~\cite{zhang2016pairwise} as a maximization problem, so \cite{zhang2016pairwise}~in fact considers block-coordinate descent as the dual is a minimization problem. However, the approach of~\cite{zhang2016pairwise} can be easily adapted to the minimization version of QAP, which we consider in our descriptions.} (BCA) in the dual problem. BCA~\cite{werner2020relative} is a well-known iterative method for (approximate) maximization of (generally constrained) multivariate functions. In each iteration, this method chooses a subset (i.e., a block) of variables and maximizes the objective over this subset of variables while keeping the other variables constant and staying within the feasible set. By repeating this iteration for different subsets, one attains a better and better solution to the original optimization problem. However, this method need not attain (or even converge to) the optimal solution of the problem.

In~\cite{zhang2016pairwise}, BCA was used to approximately optimize a dual linear programming (LP) relaxation of the QAP. The novel LP relaxation considered in~\cite{zhang2016pairwise} has the following features. If the dual linear program is restricted to one set of variables, it corresponds to the dual LP relaxation of the \emph{weighted constraint satisfaction proble}m (WCSP, equivalent to the maximum a posteriori (MAP) inference problem in graphical models~\cite{savchynskyy2019discrete,werner2007linear}). If the dual linear program is restricted to another set of variables, it becomes the LP formulation of the LAP. To obtain bounds on the optimal value of the QAP and good proposals for its solution, \cite{zhang2016pairwise}~introduced the \emph{Hungarian Belief Propagation} (Hungarian-BP) algorithm that (approximately) optimizes the dual linear program by BCA: the subproblem corresponding to the LAP is solved exactly by the Hungarian method~\cite{kuhn1955hungarian} and the dual variables corresponding to the WCSP subproblem are improved by BCA algorithm MPLP~\cite{globerson2007fixing}. After each iteration of Hungarian-BP, an assignment (i.e., a feasible solution for the QAP) is generated by solving the LAP subproblem. If the assignment is not proved to be optimal for the QAP (by comparing its cost to the bound given by the dual objective), \cite{zhang2016pairwise}~employs Hungarian-BP in a branch-and-bound scheme to either prove optimality or find a better assignment. To speed up this method for the cost of possible non-optimality of the returned assignment, one can limit the number of explored branches.

A related method was used in~\cite{hutschenreiter2021fusion} where both subproblems were optimized only approximately by BCA and, similarly to~\cite{zhang2016pairwise}, an assignment is generated after each iteration. More importantly, \cite{hutschenreiter2021fusion}~found that the quality of the generated assignments is significantly improved if their sequence is gradually fused into a single one. This method constitutes the current state of the art for computer-vision instances.

For practical purposes, especially in computer vision~\cite{haller2022comparative}, the QAP can be generalized to the setting where the sets do not have the same size and some elements may remain unassigned (possibly for some cost). This gives rise to the \emph{incomplete QAP} (IQAP) where one seeks an injective partial mapping (also called partial bijection) from one set to another that minimizes an objective consisting of both unary and binary functions, as in the QAP. Again, if the objective of the IQAP is restricted to contain only unary functions, one obtains the polynomially solvable \emph{incomplete LAP} (ILAP). 
The IQAP has been formally defined and a reduction between it and the QAP has been shown in~\cite[Section~A1]{haller2022comparative}.
However, such a straightforward reduction (based on introducing multiple additional `dummy' labels) need not be competitive with tackling the IQAP directly if the underlying problem is large-scale and sparse. In contrast, the formulation of the optimization problem in~\cite{hutschenreiter2021fusion} in fact corresponds to the IQAP (although it is referred to as the QAP in~\cite{hutschenreiter2021fusion}). Naturally, the LP relaxation of IQAP contains a subproblem corresponding to the ILAP instead of the LAP.

Here, we combine the approaches~\cite{zhang2016pairwise,hutschenreiter2021fusion} with a recent result in the theory of BCA from~\cite{werner2020relative}. It was shown that, when solving a block-subproblem, one should choose an optimizer from the relative interior of the set of optimizers, which is called the \emph{relative-interior rule}. In~\cite{werner2020relative}, this rule is shown not to be worse (in a precise sense) than any other choice of the optimizer. However, even with the relative-interior rule, the fixed points of BCA need not be optimal.

\paragraph{Contribution}
In our work, we augment the approach of~\cite{zhang2016pairwise} by choosing an optimizer from the relative interior of the set of optimizers for the ILAP subproblem. To this end, we propose and implement an algorithm that is capable of providing such a solution. Our motivation for this is that, by the results of~\cite{werner2020relative}, the fixed points of BCA algorithm conforming to the relative-interior rule are not worse than if this rule is not adhered to. Our experimental results indicate that such an approach can indeed typically attain better dual objective (i.e., a lower bound on the optimal solution). Moreover, optimizing the ILAP subproblem exactly can improve the bound in cases where BCA cannot. Our method thus extends a long line of research focused on computing bounds on the optimal value of the QAP~\cite{adams1994improved,gilmore1962optimal,lawler1963quadratic,burkard2009assignment,hahn1998lower}.

\paragraph{Structure of the Paper}
We proceed as follows. We begin in Section~\ref{se:LAP} by defining the LAP and characterizing solutions of its (dual) LP formulation that are in the relative interior of the set of optimizers. These results allow us to devise an algorithm for computing such solutions. In particular, we constructively show that one can compute a relative-interior solution from an arbitrary optimal solution in linear time (assuming that an optimal assignment for the LAP is also available). Section~\ref{se:iLAP} presents the ILAP together with a linear-time reduction of the ILAP to the LAP. To be able to compute relative-interior solutions for the LP formulation of the ILAP, we define a closed-form mapping that computes such a solution based on a relative-interior solution of the LP formulation of the constructed LAP. Next, Section~\ref{se:iQAP} formally defines the IQAP and its considered LP relaxation whose subproblems correspond to LP relaxation of the WCSP and LP formulation of the ILAP. Finally, Section~\ref{se:compared_methods} precisely describes the compared methods and overviews our experimental results.

Up to certain details, we follow the notation of~\cite{hutschenreiter2021fusion} throughout our exposition.

\section{The LAP and Relative-Interior Solution}\label{se:LAP}

Let~$\V$ be a finite set of \emph{vertices} and~$\L$ be a finite set of \emph{labels} such that~$|\V|=|\L|$ and~$\V\cap\L=\emptyset$. Next, for each~$v\in\V$, let~$\L_v\subseteq \L$ be the set of \emph{allowed labels} for vertex~$v$ and~$\theta_v\colon \L_v\rightarrow\mathbb{R}$ be a cost function. The \emph{linear assignment problem} (LAP)~\cite{burkard2009assignment} is the optimization problem
\begin{subequations}\label{eq:LAP}
\begin{align}
    \min \sum_{v\in\V} \theta_v(x_v) \span \label{eq:LAP:a}\\
    \forall v\in\V:\;& x_v\in\L_v \\
    \forall \l\in\L:\;& \sum_{v\in\V} \iverson{x_v=\l}=1\label{eq:LAP:c}
\end{align}
\end{subequations}
where~$\iverson{\cdot}$ denotes the Iverson bracket, i.e., $\iverson{\Psi}=1$ if~$\Psi$ is true and $\iverson{\Psi}=0$ otherwise. In words, the task is to find a bijection~$x\colon \V\rightarrow\L$ such that each vertex is assigned an allowed label and the objective~\eqref{eq:LAP:a} is minimized.\footnote{As usual, we will write $x\colon \V\rightarrow\L$ and $x \in \L^\V$ interchangeably in the sequel, understanding that $\L^\V$ is the set of all mappings from $\V$ to $\L$.} Note, \eqref{eq:LAP}~need not always be feasible, but we will assume in the sequel that it is.

\begin{remark}
Frequently~\cite{kuhn1955hungarian}, we have~$\L_v=\L$ for each~$v\in\V$ and the problem is thus always feasible. Although one can assume this without loss of generality (by setting a high cost~$\theta_v(\l)$ to disallowed labels~\cite{burkard2009assignment}), we use the formalism above to be consistent with recent literature on the QAP~\cite{hutschenreiter2021fusion} that is introduced later.
\end{remark}

The LAP has a natural LP formulation that can be stated as the left-hand problem of the primal-dual pair
\begin{subequations}\label{eq:LAP_LP_relax}
\begin{align}
    && \min \; \sum_{\substack{v\in\V\\\l\in\L_v}} \theta_v(\l)\mu_v(\l) \span  & \max \; \sum_{v\in\V} \alpha_v+\sum_{\l\in\L} \beta_\l \span \\
    \forall v\in\V, \,\l\in\L_v:\;&& \mu_v(\l)&\geq 0 & \alpha_v+\beta_\l&\leq\theta_v(\l)\label{eq:LAP_LP_relax:b}\\
    \forall v\in\V:\;&& \sum_{\l\in\L_v} \mu_v(\l)&=1 & \alpha_v&\in\mathbb{R}\label{eq:LAP_LP_relax:c}\\
    \forall \l\in\L:\;&& \sum_{v\in\V_\l} \mu_v(\l)&= 1 & \beta_\l&\in\mathbb{R}\label{eq:LAP_LP_relax:d}
\end{align}
\end{subequations}
where~$\V_\l=\{v\in\V\mid\l\in\L_v\}$ is the set of vertices for which the label~$\l$ is allowed. On the right, we wrote the dual linear program. We note that the dual variables are written on the same lines as the primal constraints to which they correspond and vice versa -- this applies to each primal-dual pair that we consider in this paper.

It is known~\cite{birkhoff1946tres,papadimitriou1998combinatorial,burkard2009assignment} that the constraint matrix of the primal~\eqref{eq:LAP_LP_relax} (on the left) is totally unimodular, so the vertices of the polyhedron that constitutes the feasible set are integral and correspond to feasible solutions of~\eqref{eq:LAP}. Since the minimum of a linear function on a polyhedron is always attained in at least one vertex of the polyhedron, the optimal objectives of~\eqref{eq:LAP} and~\eqref{eq:LAP_LP_relax} coincide.

\subsection{Relative Interior of the Set of Optimal Solutions}

We will assume that the sets $\V$, $\L$, $\{\L_v\}_{v\in\V}$, and costs $\theta$ are fixed in this and the following subsection to simplify formulations of statements. In this subsection, we provide characterizations of solutions of the primal and dual~\eqref{eq:LAP_LP_relax} that are in the relative interior of the set of primal and dual optimizers, respectively.

Formally, the \emph{relative interior} of a convex set~$S$, denoted~$\ri S$, is the topological interior of~$S$ relative to the affine hull of~$S$~\cite{hiriart2004fundamentals}. In the sequel, we will not use the definition of relative interior directly, but only rely on its properties. It is important to note that, for any convex set~$S$, $\ri S\subseteq S$. Moreover, it is known that $\ri S=\emptyset$ if and only if~$S=\emptyset$~\cite{hiriart2004fundamentals}.




\begin{example}[Interior vs. relative interior]
Consider three distinct points in 3-dimensional space that do not lie on the same line. Their convex hull, a two-dimensional triangle, has empty interior as all of its points are bounding it from the surrounding space. However, its relative interior is non-empty as it consists of all its points except for the points lying on the three edges.

In general, the relative interior of a non-empty polyhedron defined by the linear inequality system $Ax\leq b$ is the set of points $x$ within this polyhedron for which the number of active inequalities (i.e., inequalities satisfied with equality) among $Ax\leq b$ is the least.
We refer the interested reader to \cite{hiriart2004fundamentals} for more details.
\end{example}

Since we assumed that~\eqref{eq:LAP} is feasible, the primal~\eqref{eq:LAP_LP_relax} is also feasible and bounded. By strong duality, the dual~\eqref{eq:LAP_LP_relax} is also feasible and bounded. Both the primal and the dual thus have a non-empty set of optimizers which implies that the relative interiors of the sets of optimizers are non-empty too.

Recall that the solutions that lie in the relative interior of optimizers of any primal-dual pair can be characterized using the \emph{strict complementary slackness} condition~\cite{zhang1994strictly}. We formulate this condition for the case of the previously stated primal-dual pair~\eqref{eq:LAP_LP_relax} in the following theorem.

\begin{theorem}[\cite{zhang1994strictly}]\label{th:strict_compl_slackness}
Let~$\mu$ and~$(\alpha,\beta)$ be feasible for the primal and dual~\eqref{eq:LAP_LP_relax}, respectively. The following are equivalent:
\begin{enumerate}
    \item[(a)] $\mu$ and~$(\alpha,\beta)$ is in the relative interior of optimal solutions of the primal and dual, respectively,
    \item[(b)] $\forall v\in\V, \,\l\in\L_v: (\mu_v(\l)>0 \iff \alpha_v+\beta_\l=\theta_v(\l))$.
\end{enumerate}
\end{theorem}

Next, we provide a more tangible description of optimal solutions from the relative interior using the notion of minimally-assignable pairs~$\{v,\l\}$. For~$v\in\V$ and~$\l\in\L_v$, we say that the pair~$\{v,\l\}$ is \emph{minimally assignable} if there exists an assignment~$x$ optimal for~\eqref{eq:LAP} with~$x_v=\l$.

\begin{proposition}[cf. Property~4 in \cite{claus2020analysis}]\label{pr:rel_int_primal}
Let~$\mu$ be feasible for the primal~\eqref{eq:LAP_LP_relax}. The following are equivalent:
\begin{enumerate}
    \item[(a)] $\mu$~is in the relative interior of optimizers of the primal~\eqref{eq:LAP_LP_relax},
    \item[(b)] $\forall v\in\V,\,\l\in\L_v:$ ($\mu_v(\l)>0$ if and only if the pair~$\{v,\l\}$ is minimally assignable).
\end{enumerate}
\end{proposition}
\begin{proof}
Let $X=\{x^1, \ldots, x^n\}\subseteq \L^\V$ be the set of assignments optimal for~\eqref{eq:LAP}. Next, define~$\mu^i$ for each~$i\in\{1, \ldots, n\}$ by
\begin{equation}\label{eq:mu_from_x}
\forall v\in\V,\,\l\in\L_v: \mu_v^i(\l)=\iverson{x^i_v=\l}
\end{equation}
and let~$M = \{\mu^1, \ldots, \mu^n\}$. Due to integrality of the primal~\eqref{eq:LAP_LP_relax}, the set of solutions optimal for the primal~\eqref{eq:LAP_LP_relax} is the convex hull of~$M$, i.e.,
\begin{equation}\label{eq:conv_M}
    \conv M = \Big\{\sum_{i=1}^na_i\mu^i\,\Big\vert\, a\in\mathbb{R}^n, \,\sum_{i=1}^n a_i=1,\, \forall i\in\{1, \ldots n\}: a_i\geq0\Big\}.
\end{equation}
By~\cite[Remark~2.1.4]{hiriart2004fundamentals}, we have that the relative interior of the polyhedron~\eqref{eq:conv_M} (i.e., the relative interior of optimal solutions of the primal~\eqref{eq:LAP_LP_relax}) is
\begin{equation}\label{eq:ri_conv_M}
    \ri \conv M = \Big\{\sum_{i=1}^na_i\mu^i\,\Big\vert\, a\in\mathbb{R}^n,\, \sum_{i=1}^n a_i=1, \,\forall i\in\{1, \ldots n\}: a_i>0\Big\}
\end{equation}
where the only difference to~\eqref{eq:conv_M} is that each coefficient~$a_i$ is required to be positive.

We begin with (a)$\implies$(b), so let~$\mu^*\in\ri\conv M$, $v\in\V$, and~$\l\in\L_v$ be arbitrary. We distinguish two cases:
\begin{itemize}
    \item If the pair~$\{v,\l\}$ is not minimally assignable, then we have~$x_v\neq\l$ for all assignments~$x$ optimal for~\eqref{eq:LAP}. Consequently, $\mu_v^i(\l)=0$ for all~$i\in\{1,\ldots,n\}$ and~$\mu_v^*(\l)=\sum_{i=1}^n a_i\mu_v^i(\l)=0$. 
    \item If the pair~$\{v,\l\}$ is minimally assignable, then there exists an assignment~$x^{i^*}$ optimal for~\eqref{eq:LAP} with~$x_v^{i^*}=\l$ and thus $\mu_v^{i^*}(\l)=1$. We have that~$\mu^*_v(\l) = \sum_{i=1}^n a_i\mu_v^i(\l)>0$ because all the terms~$a_i\mu_v^i(\l)$ are non-negative and $a_{i^*}\mu_v^{i^*}(\l)>0$ due to $a_{i^*}>0$.
\end{itemize}
We continue to prove (b)$\implies$(a). Let~$\mu'$ and~$(\alpha',\beta')$ be from the relative interior of optimizers of the primal and dual, respectively, and let~$\mu^*$ be feasible for the primal~\eqref{eq:LAP_LP_relax}. If~$\mu^*$ satisfies condition~(b), then, for all~$v\in\V$ and~$\l\in\L_v$,
\begin{equation}\label{eq:equivalences}
    \mu^*_v(\l)>0 \iff\{v,\l\} \text{ is minimally assignable} \iff \mu'_v(\l)>0 \iff \alpha'_v+\beta'_\l=\theta_v(\l)
\end{equation}
where the first equivalence is statement~(b), the second equivalence follows from (a)$\implies$(b) and the fact that~$\mu'$ is in the relative interior of optimizers of the primal. The third equivalence in~\eqref{eq:equivalences} is strict complementary slackness for~$(\alpha',\beta')$ and~$\mu'$ and follows from Theorem~\ref{th:strict_compl_slackness}. Consequently, $\mu^*$ is in the relative interior of optimizers of the primal because it satisfies strict complementary slackness with~$(\alpha',\beta')$ by composing the equivalences in~\eqref{eq:equivalences}.
\end{proof}

Combining Proposition~\ref{pr:rel_int_primal} with Theorem~\ref{th:strict_compl_slackness} results in the following corollary:

\begin{corollary}\label{co:dual}
Let~$(\alpha,\beta)$ be feasible for the dual~\eqref{eq:LAP_LP_relax}. The following are equivalent:
\begin{enumerate}
    \item[(a)] $(\alpha,\beta)$ is in the relative interior of optimizers of the dual~\eqref{eq:LAP_LP_relax},
    \item[(b)] $\forall v\in\V,\,\l\in\L_v:$ ($\alpha_v+\beta_\l=\theta_v(\l)$ if and only if the pair~$\{v,\l\}$ is minimally assignable).
\end{enumerate}
\end{corollary}

For the purposes of the sequel, we introduce the notion of perfectly matchable edges. Formally, let~$\E\subseteq\{\{v,\l\}\mid v\in\V,\,\l\in\L_v\}$ so that~$(\V\cup\L,\E)$ is a bipartite graph. By a \emph{perfect matching} in~$(\V\cup\L,\E)$, we mean a bijection~$x\colon\V\rightarrow\L$ such that~$\{v,x_v\}\in\E$ for all~$v\in\V$. An edge~$\{v,\l\}\in\E$ is \emph{perfectly matchable} in~$(\V\cup\L,\E)$ if~$(\V\cup\L,\E)$ has a perfect matching~$x$ with~$x_v=\l$.

With different choices of the edge set~$\E$, there are several connections between perfect matchings in~$(\V\cup\L,\E)$ and assignments feasible or optimal for~\eqref{eq:LAP} which we outline in the remaining part of this subsection.

For example, with~$\E=\{\{v,\l\}\mid v\in\V,\,\l\in\L_v\}$, the set of perfect matchings in~$(\V\cup\L,\E)$ coincides with the feasible set of~\eqref{eq:LAP}. With this choice of~$\E$, an edge~$\{v,\l\}\in\E$ is perfectly matchable in~$(\V\cup\L,\E)$ if and only if there exists an assignment~$x$ feasible for~\eqref{eq:LAP} with~$x_v=\l$.

Recall that, for~$(\alpha,\beta)$ feasible for the dual~\eqref{eq:LAP_LP_relax}, the bipartite graph~$(\V\cup\L,\E(\alpha,\beta))$ where
\begin{equation}\label{eq:equality_subgraph}
\E(\alpha,\beta)=\{\{v,\l\}\mid v\in\V,\,\l\in\L_v,\,\alpha_v+\beta_\l=\theta_v(\l)\}
\end{equation}
is called the \emph{equality subgraph}~\cite{mills2007dynamic,akgul1992linear}. The following lemma connects perfectly matchable and minimally-assignable edges in the context of the LAP and the equality subgraph.
\begin{lemma}\label{le:min_assignable_iff_perf_matchable}
Let~$(\alpha,\beta)$ be optimal for the dual~\eqref{eq:LAP_LP_relax}.
For all~$v\in\V$ and~$\l\in\L_v$, $\{v,\l\}$~is minimally assignable if and only if $\{v,\l\}$~is perfectly matchable in~$(\V\cup\L,\E(\alpha,\beta))$.
\end{lemma}
\begin{proof}
Let~$v\in\V$ and~$\l\in\L_v$ be arbitrary. If $\{v,\l\}$~is minimally assignable, there exists an assignment~$x$ optimal for~\eqref{eq:LAP} such that~$x_v=\l$. Analogously to~\eqref{eq:mu_from_x}, $\mu$~defined by $\mu_v(\l)=\iverson{x_v=l}$ is optimal for the primal~\eqref{eq:LAP_LP_relax} and thus satisfies complementary slackness with the dual-optimal solution~$(\alpha,\beta)$. Consequently, $\alpha_{v'}+\beta_{x_{v'}}=\theta_{v'}(x_{v'})$ holds for all~$v'\in\V$, hence~$\{v',x_{v'}\}\in\E(\alpha,\beta)$ for all~$v'\in\V$, $x$~is a perfect matching in the equality subgraph, and $\{v,\l\}$~is perfectly matchable.

For the converse relation, if $\{v,\l\}$~is perfectly matchable, there exists a perfect matching~$x$ in~$(\V\cup\L,\E(\alpha,\beta))$ with~$x_v=\l$ by definition. This assignment~$x$ is optimal for~\eqref{eq:LAP} since~$\mu$ defined as above satisfies complementary slackness with~$(\alpha,\beta)$.
\end{proof}

We are now able to formulate Theorem~\ref{th:summary}, which is the important result of this subsection.

\begin{theorem}\label{th:summary}
Let~$(\alpha,\beta)$ be feasible for the dual~\eqref{eq:LAP_LP_relax}.
\begin{enumerate}[label={(\alph*)},ref={\thetheorem\alph*}]
    \item $(\alpha,\beta)$ is optimal for the dual if and only if the equality subgraph~$(\V\cup\L,\E(\alpha,\beta))$ has a perfect matching,
    \item $(\alpha,\beta)$~is in the relative interior of optimizers of the dual if and only if the equality subgraph~$(\V\cup\L,\E(\alpha,\beta))$ has a perfect matching and each edge of the equality subgraph is perfectly matchable.\label{th:summary:b}
\end{enumerate}
\end{theorem}
\begin{proof}
Statement~(a) follows from complementary slackness: for a perfect matching~$x$, we define a feasible solution~$\mu$ of primal~\eqref{eq:LAP_LP_relax} as in~\eqref{eq:mu_from_x}. This pair of solutions satisfies complementary slackness, so $\mu$~is optimal for the primal and~$(\alpha,\beta)$ is optimal for the dual. Conversely, if~$(\alpha,\beta)$ is optimal for the dual and~$x$ is an optimal assignment for~\eqref{eq:LAP}, then~$\mu$ defined by~\eqref{eq:mu_from_x} based on~$x$ is optimal for the primal and satisfies complementary slackness conditions, i.e., $\alpha_v+\beta_{x_v}=\theta_v(x_v)$ holds for all~$v\in\V$ due to~$\mu_v(x_v)=1>0$, so~$x$ is a perfect matching in~$(\V\cup\L,\E(\alpha,\beta))$.

Statement~(b) is obtained by combining Lemma~\ref{le:min_assignable_iff_perf_matchable} with Corollary~\ref{co:dual}.
\end{proof}

\subsection{Obtaining a Relative-Interior Solution from an Optimal Solution}

Let us now focus on obtaining a solution belonging to the relative interior of optimal solutions of the dual~\eqref{eq:LAP_LP_relax}. For this, we assume that an arbitrary dual-optimal solution~$(\alpha,\beta)$ and an assignment~$x\colon \V\rightarrow \L$ optimal for~\eqref{eq:LAP} are at our disposal. Both can be obtained, e.g., using the Hungarian method for solving the LAP~\cite{kuhn1955hungarian}. Our method for computing a relative-interior solution is based on changing the given dual solution~$(\alpha,\beta)$ so that it remains optimal but non-perfectly-matchable edges are removed from the equality subgraph by making the corresponding dual constraint~\eqref{eq:LAP_LP_relax:b} hold with strict inequality. In other words, only perfectly matchable edges remain in the equality subgraph, which corresponds to being in the relative interior of dual optimizers by Theorem~\ref{th:summary:b}. Next, we focus on how to perform this task in detail.

We base our method for obtaining a relative-interior solution on the paper~\cite{tassa2012finding} where (aside from other results) a method for finding all perfectly matchable edges in a bipartite graph was proposed. Although the technique of~\cite{tassa2012finding} is applicable to any bipartite graph~${(\V\cup\L,\E)}$, we will apply it to the equality subgraph, i.e., we have~$\E=\E(\alpha,\beta)$.

Following~\cite{tassa2012finding}, the first step of the procedure is to construct the directed graph~$(\V,\F_x(\E))$ where\footnote{As usual, we denote edges of undirected graphs as (unordered) 2-element sets (e.g.,~$\{u,v\}$) and edges of directed graphs as (ordered) 2-tuples (e.g.,~$(u,v)$).}
\begin{equation}\label{eq:directed_graph_for_perfectly_matchable}
    \F_x(\E) = \{(u,v)\mid u,v\in\V,\,u\neq v,\,\{u,x_v\}\in \E\}.
\end{equation}
The motivation for constructing the graph~$(\V,\F_x(\E))$ is the following. An edge~$\{v,\l\}$ is perfectly matchable in~$(\V\cup\L,\E)$ if and only if~$x_v=\l$ or there exists an alternating cycle in~$(\V\cup\L,\E)$ w.r.t.~$x$ that contains~$\{v,\l\}$~\cite{tassa2012finding}.\footnote{Recall~\cite{tassa2012finding} that an alternating cycle in~$(\V\cup\L,\,\E)$ w.r.t.~$x$ is a sequence of (non-repeating) edges~$\{v_1,\l_1\}, \{v_2,\l_2\},\ldots, \{v_{2n},\l_{2n}\}$ from~$\E$ satisfying the following three conditions: (i) for all odd~$i$: $\l_{i}=\l_{i+1}$ and~$x_{v_i}\neq \l_i$, (ii) for all even~$i$: $v_{i}=v_{i+1}$ and~$x_{v_i}=\l_i$, (iii) $v_{2n}=v_1$.} Each alternating cycle w.r.t.~$x$ in the undirected graph~$(\V\cup\L,\E)$ corresponds to a directed cycle in~$(\V,\F_x(\E))$ and vice versa. We thus have the following result from~\cite{tassa2012finding}.

\begin{theorem}[\cite{tassa2012finding}]\label{th:tassa}
Let~$(\V\cup\L,\E)$ be a bipartite graph with a perfect matching~$x$. Let~$\{v,\l\}\in\E$ and~$u\in\V$ be such that~$x_u=\l$. Edge~$\{v,\l\}$ is perfectly matchable in~$(\V\cup\L,\E)$ if and only if~$x_v=\l$ or~$(v,u)\in\F_x(\E)$ is a part of a directed cycle in~$(\V,\F_x(\E))$.
\end{theorem}

To simplify notation in this subsection, we abbreviate~$\F=\F_x(\E)$ if there is no ambiguity. First, we focus on deciding which edges of~$(\V,\F)$ belong to some directed cycle. As proposed in~\cite{tassa2012finding}, this can be achieved by computing the strongly connected components of~$(\V,\F)$, which can be done in linear time~\cite{tarjan1972depth,sharir1981strong}. We remind\footnote{For a more detailed overview of graph-theoretical notions, we refer to the book~\cite{thulasiraman1992graphs} or papers~\cite{tassa2012finding,tarjan1972depth}.} the reader that a directed graph is \emph{strongly connected} if there exists a directed path between each ordered pair of its vertices. The \emph{strongly connected components} of a directed graph are its maximal subgraphs (w.r.t.\ their set of vertices) that are strongly connected. Since the strongly connected components are vertex-induced subgraphs, we will identify the strongly connected components with the sets of vertices that induce them. Let~$\{V_1, \ldots, V_n\}$ be the strongly connected components of~$(\V,\F)$ so that $\{V_1, \ldots, V_n\}$~is a partition of~$\V$. The \emph{condensation} of~$(\V,\F)$ is a directed acyclic graph with the set of vertices~$\{V_1, \ldots, V_n\}$ (denoted by subsets of vertices of the original graph). The condensation contains an edge~$(V_i,V_j)$ if and only if $i\neq j$ and ${\exists u\in V_i},\,v\in V_j: (u,v)\in\F$. Since the condensation is acyclic, it has a topological ordering, i.e., an ordering in which every edge leads from an earlier to a later vertex in the ordering \cite[Section~5.7]{thulasiraman1992graphs}.

Let us order the strongly connected components such that~$\{V_1, \ldots, V_n\}$ is a topological ordering of the condensation. To obtain a relative-interior solution (based on the previously mentioned optimal solution~$(\alpha,\beta)$), we process the components in a reversed topological order (i.e., an order in which each component is processed before any of its predecessors in the condensation) and, sequentially for each~$V_i$ with ingoing edges in the condensation, update
\begin{subequations}\label{eq:update_dual_variables}
\begin{align}
    \forall v \in V_i:&\; \alpha_v := \alpha_v +\delta/2\\
    \forall \l \in L_i:&\; \mathrlap{\beta_\l}\phantom{\alpha_v} := \mathrlap{\beta_\l}\phantom{\alpha_v} -\delta/2
\end{align}
\end{subequations}
where\footnote{Notice that~\eqref{eq:update_dual_variables} resembles the update of dual variables in the Hungarian method~\cite[Page~94]{kuhn1955hungarian}, but there are clear differences. E.g., our updates do not change the dual objective and make some of the dual constraints inactive (i.e., hold with strict inequality instead of equality) whereas in the Hungarian method, the updates improve dual objective and make some dual constraints active.} $L_i = \{x_v\mid v \in V_i\}$ and
\begin{equation}\label{eq:delta}
    \delta = \begin{cases}
    \min\{ \theta_v(\l)-\alpha_v-\beta_\l \mid v \in V_i,\, \l\in \L_v\setminus L_i \} & \text{ if this set is non-empty}\\
    1 & \text{ otherwise }
    \end{cases}.
\end{equation}

We provide a formal overview of this procedure in Algorithm~\ref{al:ri-solution}. Correctness and time complexity of this algorithm is given by Theorem~\ref{th:correctness_inductive}, which requires an auxiliary lemma.

\begin{algorithm}[t]\SetAlgoVlined
\SetKwInOut{Output}{output} 
\Input{instance of the LAP, optimal solution~$(\alpha,\beta)$ for dual~\eqref{eq:LAP_LP_relax}, assignment~$x$ optimal for~\eqref{eq:LAP}}
\Output{solution from the relative interior of optimizers of the dual~\eqref{eq:LAP_LP_relax}}
Construct graph~$(\V,\F)$ by computing~$\F:=\F_x(\E(\alpha,\beta))$ (see~\eqref{eq:directed_graph_for_perfectly_matchable}).\\
Compute the strongly connected components and the condensation of~$(\V,\F)$.\\
Find a topological ordering~$\{V_1, \ldots, V_n\}$ of the condensation of~$(\V,\F)$.\\
\For{$i\in\{n, \ldots, 1\}$ \normalfont{(in decreasing order)}}{
\If{\normalfont{component~$i$ has ingoing edges in the condensation}}{
Perform updates~\eqref{eq:update_dual_variables} where~$\delta$ is~\eqref{eq:delta} and~$L_i = \{x_v\mid v \in V_i\}$.
}
}
\Return{$(\alpha,\beta)$}
\caption{Computing a solution from the relative-interior of optimizers of dual~\eqref{eq:LAP_LP_relax}.}
\label{al:ri-solution}
\end{algorithm}

\begin{lemma}\label{le:correctness}
Let~$(\alpha,\beta)$ be optimal for the dual~\eqref{eq:LAP_LP_relax}. Let~$\{V_1,\ldots,V_n\}$ be the strongly connected components of~$(\V,\F_x(\E(\alpha,\beta))$. Let~$i\in\{1,\ldots,n\}$. Perform the update on lines 5--6 of Algorithm~\ref{al:ri-solution} for~$i$ and denote the resulting values of the dual variables by~$(\alpha',\beta')$. If there are no outgoing edges from the component~$V_i$ in the condensation of~$(\V,\F_x(\E(\alpha,\beta)))$, then
\begin{enumerate}[label={(\alph*)},ref={\thetheorem\alph*}]
    \item $(\alpha',\beta')$ is dual optimal,\label{le:correctness:a}
    \item $\E(\alpha',\beta') = \E(\alpha,\beta) \setminus \{\{v,\l\}\in \E(\alpha,\beta)\mid \l \in L_i,\,v\notin V_i \}$,\label{le:correctness:b}
    \item $\F_x(\E(\alpha',\beta'))=\F_x(\E(\alpha,\beta))\setminus \{(v,u)\in\F_x(\E(\alpha,\beta))\mid v\in \V\setminus V_i,\,u\in V_i\}$.\footnote{Condition (c) is implied by (b) following the definition~\eqref{eq:directed_graph_for_perfectly_matchable}. We use (c) for the sake of the proof though.
    }\label{le:correctness:c}
\end{enumerate}
\end{lemma}
\begin{proof}
If condition on line~5 of Algorithm~\ref{al:ri-solution} is not satisfied, then the strongly connected component~$V_i$ is an isolated vertex in the condensation because it has no outgoing edges (by our assumption in the lemma) and no ingoing edges (by condition on line~5). Statements (a)-(c) are thus trivially satisfied by~$(\alpha',\beta')=(\alpha,\beta)$.

For the remaining part, let condition on line~5 be satisfied, i.e., $(u,v)\in\F_x(\E(\alpha,\beta))$ for some~$u\in \V\setminus V_i$ and~$v\in V_i$. Consequently, we have~$\{u,x_v\}\in\E(\alpha,\beta)$, i.e., there is at least one edge between~$\V\setminus V_i$ and~$L_i$ in $\E(\alpha,\beta)$.

Next, see that there are no edges between~$V_i$ and~$\L\setminus L_i$ in~$\E(\alpha,\beta)$ -- for contradiction, let~$v\in V_i$, $\l\in\L\setminus L_i$, and~$\{v,\l\}\in\E(\alpha,\beta)$. Denoting by~$u\in \V$ the vertex with~$x_u=\l$ yields that~$(v,u)\in\F_x(\E(\alpha,\beta))$ by~\eqref{eq:directed_graph_for_perfectly_matchable}. See that~$u\in\V\setminus V_i$ because $\l\in\L\setminus L_i$, which implies that there is an outgoing edge from the component~$V_i$. This is contradictory with the fact that~$V_i$ has no outgoing edges by our assumption in this lemma. This also implies that~$\delta > 0$ because~$\alpha_v+\beta_\l<\theta_v(\l)$ holds for each~$v\in V_i$ and~$\l\in\L_v \setminus L_i$.

We proceed to show by case analysis that the new values~$(\alpha',\beta')$ satisfy statements~(b) and~(c). Let~$v\in\V$ and~$\l\in\L_v$.
\begin{itemize}
    \item If~$v\in V_i$ and~$\l\in\L_v\setminus L_i$, then~$\{v,\l\}\notin\E(\alpha,\beta)$, as discussed previously. By definition of~$\delta$, we have $\delta\leq \theta_v(\l)-\alpha_v-\beta_\l$ and thus~$\alpha_v'+\beta_\l'=\alpha_v+\beta_\l+\delta/2<\theta_v(\l)$, i.e., $\{v,\l\}\notin\E(\alpha',\beta')$.
    \item If~$v\in V_i$ and~$\l\in L_i$ (or $v\in\V\setminus V_i$ and~$\l\in\L \setminus L_i$), then~$\{v,\l\}\in\E(\alpha,\beta)$ if and only if~$\{v,\l\}\in\E(\alpha',\beta')$ due to~$\alpha_v+\beta_\l=\alpha_v'+\beta_\l'$ by~\eqref{eq:update_dual_variables}.
    \item If~$v\in\V \setminus V_i$ and~$\l\in L_i$, then the vertex~$u\in V_i$ with~$x_u=\l$ is in a different component than~$v$. By~$\delta>0$, we have that~$\alpha_v'+\beta_\l'<\alpha_v+\beta_\l\leq\theta_v(\l)$, so~$\{v,\l\}\notin \E(\alpha',\beta')$ and~$(v,u)\notin\F_x(\E(\alpha',\beta'))$.
\end{itemize}

Optimality of~$(\alpha',\beta')$ follows from the fact that the objective does not change by the update~\eqref{eq:update_dual_variables}, i.e., $\sum_{v\in\V}\alpha_v+\sum_{\l\in\L}\beta_\l=\sum_{v\in\V}\alpha_v'+\sum_{\l\in\L}\beta_\l'$. Feasibility follows from the previous case analysis.
\end{proof}

\begin{theorem}\label{th:correctness_inductive}
Algorithm~\ref{al:ri-solution} returns a dual-optimal solution from the relative interior of optimizers. Moreover, the time complexity of Algorithm~\ref{al:ri-solution} is~$O(\sum_{v\in\V}|\L_v|)$, i.e., linear in the size of the input.
\end{theorem}
\begin{proof}
To show correctness, we proceed by induction based on Lemma~\ref{le:correctness}. For this, let us denote by~$(\alpha^{n+1},\beta^{n+1})$ the initial values of the dual variables and by~$(\alpha^i,\beta^i)$ their values after the update on line~6 was performed with~$i$. In case that the update is skipped for~$i$ (due to unsatisfied condition on line~5), we define~$(\alpha^i,\beta^i)=(\alpha^{i+1},\beta^{i+1})$. We will prove that the following holds for any~$i\in\{1,\ldots,n\}$:
\begin{enumerate}
    \item[(a)] $(\alpha^i,\beta^i)$ is dual optimal,
    \item[(b)] $\displaystyle\E^i=\E^{n+1}\setminus \bigcup\limits_{j=i}^n\{\{v,\l\}\in \E^{n+1}\mid \l\in L_j, v\notin V_j\}$,
    \item[(c)] $\displaystyle\F_x(\E^i)=\F_x(\E^{n+1})\setminus\bigcup\limits_{j=i}^n\{(v,u)\in \F_x(\E^{n+1})\mid v\in \V\setminus V_j,\,u\in V_j \}$
\end{enumerate}
where we abbreviated~$\E(\alpha^i,\beta^i)$ to~$\E^i$ for each~$i\in\{1, \ldots, n+1\}$.

First, we show that the assumptions of Lemma~\ref{le:correctness} are satisfied whenever we process any component~$V_i$. For the base case, note that there can be no outgoing edges from~$V_n$ as it is last in the topological ordering, so the assumptions of Lemma~\ref{le:correctness} are satisfied. For the inductive step, based on statement~(c) in Lemma~\ref{le:correctness}, the edges leading to processed components are removed, so whenever a component~$V_i$ should be processed, it has no outgoing edges since all its successor components must have been processed beforehand.

The statements (a)-(c) in this proof thus follow from inductively combining the statements (a)-(c) in Lemma~\ref{le:correctness}. When all components are processed, we have that~$(\alpha^1,\beta^1)$ is optimal by~(a) and
\begin{subequations}
\begin{align}
    \E^1&=\{\{v,\l\}\in \E^{n+1}\mid \exists j\in\{1,\ldots,n\}:(v,\l)\in V_j\times L_j\}\\
    \F_x(\E^1)&=\{(u,v)\in \F_x(\E^{n+1})\mid \exists j\in\{1,\ldots,n\}: u,v\in V_j \}\label{eq:th_f_n}
\end{align}
\end{subequations}
by~(b) and~(c), respectively. Note that~$\F_x(\E^{1})\subseteq \F_x(\E^{n+1})$, so~$(\V,\F_x(\E^{1}))$ is a subgraph of the directed graph~$(\V,\F_x(\E^{n+1}))$. Moreover, $\F_x(\E^{1})$~contains precisely those edges that are within the strongly connected components of~$(\V,\F_x(\E^{n+1}))$. By Theorem~\ref{th:tassa}, each edge of~$(\V,\F_x(\E^1))$ belongs to some directed cycle and thus each edge in~$\E^1=\E(\alpha^1,\beta^1)$ is perfectly matchable. By Theorem~\ref{th:summary:b}, $(\alpha^1,\beta^1)$ is in the relative interior of optimizers of the dual.

Concerning the time complexity, this is clear for the construction of~$(\V,\F_x(\E^{n+1}))$ and also for lines~2--3 that can be performed in linear time, e.g., by Tarjan's algorithm~\cite{tarjan1972depth}. Next, for each~$i$, $\delta$ can be computed in at most~$O(\sum_{v\in V_i}|\L_v|)$ operations and update of~$(\alpha,\beta)$ can be performed in~$O(|V_i|)$ operations. All in all, the loop on lines~4--6 takes at most~$O(\sum_{v\in\V}|\L_v|)$ operations. Note that~$\sum_{v\in\V}|\L_v|\geq|\V|$.
\end{proof}

\begin{remark}
Results related to what we reviewed and derived above were also discovered in the constraint programming community. For example, the problem of identifying the perfectly matchable edges is equivalent to arc-consistency filtering in the alldifferent constraint, as shown in \cite{regin1994}. In~\cite{german2017arc}, the authors utilize the dual linear program to filter domain elements in constraints to achieve arc consistency and also discuss applications to weighted CSPs using the reduced costs. The paper \cite{claus2020analysis} further extends and details these ideas for particular constraints.

We highlight that both \cite{german2017arc} and \cite{claus2020analysis} (somewhat informally) mention the notion of an `interior' solution that can be used to most efficiently filter out domain elements. This precisely corresponds to the relative interior that we consider.
\end{remark}

\begin{remark}
The reverse topological order when performing updates~\eqref{eq:update_dual_variables} is needed to preserve feasibility of the solution throughout the algorithm. If the topological order is not reversed, increasing the $\alpha$ variables and decreasing the $\beta$ variables may lead to an infeasible solution.
\end{remark}

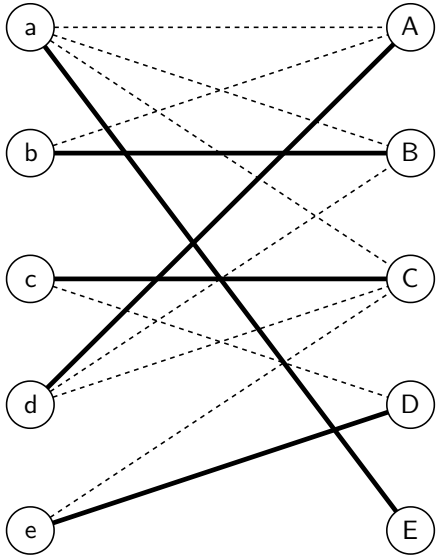
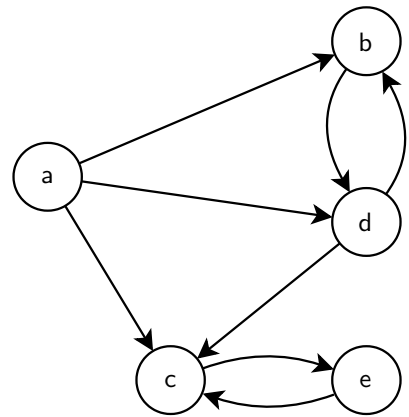
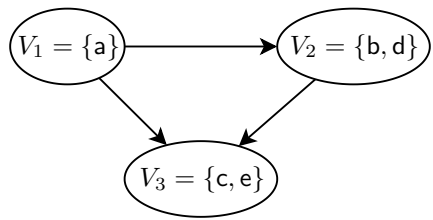
\begin{figure}[]
\centering
\begin{subfigure}[c]{\textwidth}
   \centering
     \begin{tabular}{c|c|c|c|c|c|}
        & \val{A} & \val{B} & \val{C} & \val{D} & \val{E} \\ 
        & $\beta_{\val{A}}=1$ & $\beta_{\val{B}}=1$ & $\beta_{\val{C}}=1$ & $\beta_{\val{D}}=4$ & $\beta_{\val{E}}=4$ \\ \hline
      \makecell{\val{a} \\ $\alpha_{\val{a}}=2$} &  3\cellcolor{lightgray}& 3\cellcolor{lightgray}& 3\cellcolor{lightgray}& 7 & \circled{6}\cellcolor{lightgray}\\ \hline  
      \makecell{\val{b}\\ $\alpha_{\val{b}}=2$} & \circled{3}\cellcolor{lightgray}& \circled{3}\cellcolor{lightgray}& 9 & 9 & 8 \\ \hline
      \makecell{\val{c}\\ $\alpha_{\val{c}}=3$} & 9 & 10& \circled{4}\cellcolor{lightgray}& \circled{7}\cellcolor{lightgray}& 11 \\ \hline
      \makecell{\val{d}\\ $\alpha_{\val{d}}=3$} & \circled{4}\cellcolor{lightgray}& \circled{4}\cellcolor{lightgray}& 4\cellcolor{lightgray}& 8 & 11\\ \hline
      \makecell{\val{e}\\$\alpha_{\val{e}}=3$} & 8 & 9 & \circled{4}\cellcolor{lightgray}& \circled{7}\cellcolor{lightgray}& 13 \\ \hline
     \end{tabular}
   \caption{Costs~$\theta$ for individual edges of the complete bipartite graph~$(\V\cup\L,\E)$ and initial values of dual variables~$(\alpha,\beta)$. Edges that are present in the equality subgraph are highlighted in gray. The costs of the edges that are perfectly matchable in the equality subgraph (i.e., minimally assignable for the LAP defined by these costs) are circled.}
   \label{fig:table_initial}
\end{subfigure}

\begin{subfigure}[c]{.45\textwidth}
    \centering
    \begin{tikzpicture}[baseline={(0,0)}]
        \node (pic) at (0,0) {\includegraphics[width=162pt]{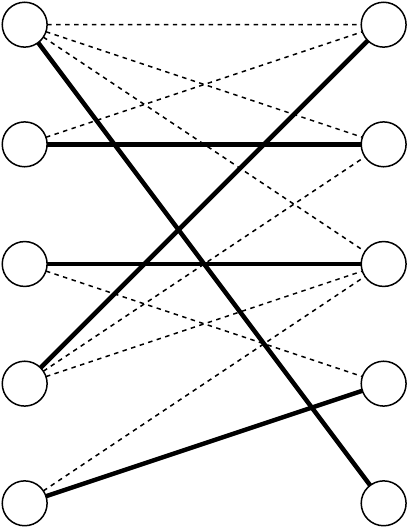}};
        \node[] at (-71pt,95pt) {\val{a}};
        \node[] at (-71pt,48.5pt) {\val{b}};
        \node[] at (-71pt,0pt) {\val{c}};
        \node[] at (-71pt,-46.5pt) {\val{d}};
        \node[] at (-71pt,-95pt) {\val{e}};
        
        \node[] at (71pt,96pt) {\val{A}};
        \node[] at (71pt,48.5pt) {\val{B}};
        \node[] at (71pt,1pt) {\val{C}};
        \node[] at (71pt,-46.5pt) {\val{D}};
        \node[] at (71pt,-94pt) {\val{E}};
    \end{tikzpicture}
    \caption{Equality subgraph for the initial values of the dual variables with a highlighted perfect matching.}
    \label{fig:graph_with_matching} 
\end{subfigure}\hspace{0.05\textwidth}
\begin{subfigure}[c]{.45\textwidth}
    \centering
    \begin{tikzpicture}[baseline={(0,0)}]
        \node (pic) at (0,0) {\includegraphics[width=150pt]{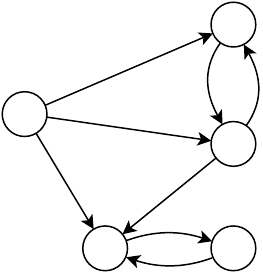}};
        \node[] at (-61pt,14pt) {\val{a}};
        \node[] at (58pt,65pt) {\val{b}};
        \node[] at (58pt,-3pt) {\val{d}};
        \node[] at (-15pt,-63pt) {\val{c}};
        \node[] at (58pt,-63pt) {\val{e}};
    \end{tikzpicture}
    \caption{Oriented graph~$(\V,\F)$ defined by the equality subgraph and perfect matching~$x$ in Figure~\ref{fig:graph_with_matching}.}
    \label{fig:oriented_graph} 
\end{subfigure}
\begin{subfigure}[c]{\textwidth}
    \centering
    \begin{tikzpicture}[baseline={(0,0)}]
        \node (pic) at (0,0) {\includegraphics[width=160pt]{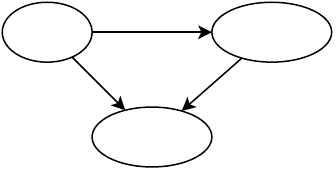}};
        \node[] at (-57pt,25.75pt) {$V_1=\{\val{a}\}$};
        \node[] at (50pt,25.75pt) {$V_2=\{\val{b},\val{d}\}$};
        \node[] at (-7pt,-24.25pt) {$V_3 = \{\val{c},\val{e}\}$};
    \end{tikzpicture}
    \caption{Condensation of the oriented graph~$(\V,\F)$ from Figure~\ref{fig:oriented_graph}.}
    \label{fig:condensation} 
\end{subfigure}
\caption{Illustrations to Example~\ref{ex:lap_ri}.}
\label{fig:example}
\end{figure}

\begin{example}\label{ex:lap_ri}
Let us consider the LAP with~$\V=\{\val{a},\val{b},\val{c},\val{d},\val{e}\}$, $\L=\{\val{A},\val{B},\val{C},\val{D},\val{E}\}$, costs~$\theta$ defined in Figure~\ref{fig:table_initial}, and a dual-optimal solution~$(\alpha,\beta)$ defined by~$\alpha=(2,2,3,3,3)$, $\beta=(1,1,1,4,4)$ (also indicated in Figure~\ref{fig:table_initial}). Here, we assume for simplicity that~$\L_v=\L$ for each~$v\in\V$. Following Algorithm~\ref{al:ri-solution}, we will show how to change this dual-optimal solution~$(\alpha,\beta)$ (which is in the relative boundary) to obtain a solution that belongs to the relative interior of dual optimizers.

The equality subgraph~$(\V\cup\L,\E(\alpha,\beta))$ with a perfect matching~$x$ is shown in Figure~\ref{fig:graph_with_matching}. Based on the equality subgraph and the chosen perfect matching, one can construct the directed graph~$(\V,\F)$ that is shown in Figure~\ref{fig:oriented_graph}. The directed graph~$(\V,\F)$ contains 2 directed cycles, each corresponding to an augmenting cycle in~$(\V\cup\L,\E(\alpha,\beta))$:
\begin{itemize}
    \item cycle $(\val{b},\val{d})$, $(\val{d},\val{b})$ corresponding to augmenting cycle~$\{\val{b},\val{B}\}$, $\{\val{B},\val{d}\}$, $\{\val{d},\val{A}\}$, $\{\val{A},\val{b}\}$,
    \item cycle $(\val{c},\val{e})$, $(\val{e},\val{c})$ corresponding to augmenting cycle~$\{\val{c},\val{C}\}$, $\{\val{C},\val{e}\}$, $\{\val{e},\val{D}\}$, $\{\val{D},\val{c}\}$.
\end{itemize}
Consequently, all of these aforementioned edges (along with edges~$\{v,x_v\}$ for all~$v\in\V$) are perfectly matchable in the equality subgraph. All the other edges in the equality subgraph (e.g., $\{\val{a},\val{A}\}$ or~$\{\val{d},\val{C}\}$) are not perfectly matchable. The perfectly matchable edges in the equality subgraph are marked by circles in Figure~\ref{fig:table_initial}.

The condensation of the directed graph~$(\V,\F)$ is shown in Figure~\ref{fig:condensation} and the partition of the vertices~$\{V_1,V_2,V_3\}$ also corresponds to a topological ordering of the condensation. We will now process the strongly connected components in reversed topological order:
\begin{enumerate}
    \item[1.] For $i=3$, we have~$V_3=\{\val{c},\val{e}\}$, $L_3 = \{\val{C},\val{D}\}$, and~$\delta = 4$. The minimal value of~$\delta$ is attained, e.g., for~$v=\val{e}$ and~$\l=\val{A}$. By~\eqref{eq:update_dual_variables}, we increase~$\alpha_{\val{c}}$ and~$\alpha_{\val{e}}$ by~2 and decrease~$\beta_{\val{C}}$ and~$\beta_{\val{D}}$ by~2. This results in~$\alpha=(2,2,5,3,5)$ and~$\beta=(1,1,-1,2,4)$.
    \item[2.] For $i=2$, we have~$V_2=\{\val{b},\val{d}\}$, $L_3 = \{\val{B},\val{A}\}$, and~$\delta = 2$. The minimal value of~$\delta$ is attained, e.g., for~$v=\val{b}$ and~$\l=\val{E}$. Note that we use the previously updated values of dual variables to compute~$\delta$, not the initial ones. By~\eqref{eq:update_dual_variables}, we increase~$\alpha_{\val{b}}$ and~$\alpha_{\val{d}}$ by~1 and decrease~$\beta_{\val{B}}$ and~$\beta_{\val{A}}$ by~1. This yields~$\alpha=(2,3,5,4,5)$ and~$\beta=(0,0,-1,2,4)$.
    \item[3.] For~$i=1$, the strongly connected component~$V_1=\{\val{a}\}$ has no ingoing edges and is thus not processed.
\end{enumerate}
By Theorem~\ref{th:summary:b}, the current dual solution~$\alpha=(2,3,5,4,5)$, $\beta=(0,0,-1,2,4)$ is in the relative interior of dual optimizers because the equality subgraph~$(\V\cup\L,\E(\alpha,\beta))$ (computed for these values of dual variables) contains only perfectly matchable edges.
\end{example}


\section{The Incomplete LAP: Reduction to the LAP and Relative-Interior Solution}\label{se:iLAP}

In this section, we introduce the \emph{incomplete LAP} (ILAP) and a reduction from the ILAP to the LAP that makes it possible to solve the ILAP using widely available LAP algorithms. Moreover, this reduction allows us to easily construct a solution of the LP formulation of the ILAP based on a solution of the LP formulation of the LAP. Importantly, we discuss that one can easily map a relative-interior solution of the LP formulation of the LAP to a relative-interior solution of the LP formulation of the ILAP.

The ILAP is the extension of the LAP where each vertex from~$\V$ need not be assigned a label, i.e., some vertices may be left unassigned for a vertex-specific cost. As in~\cite{hutschenreiter2021fusion}, we capture this by extending the set of labels by the \emph{dummy label}~$\#$ that can be assigned to arbitrarily many vertices.

Formally, let~$\V$ be a finite set of vertices and~$\L$ be a finite set of labels with~$\#\in\L$ and~$\V\cap\L=\emptyset$. For each~$v\in\V$, $\L_v\subseteq \L$ with~$\#\in\L_v$ is the set of allowed labels for vertex~$v$. Finally, as in the LAP,~$\theta_v\colon \L_v\rightarrow\mathbb{R}$ is a cost function for each~$v\in\V$. In this setting, the ILAP is the optimization problem
\begin{subequations}\label{eq:iLAP}
\begin{align}
    \min \sum_{v\in\V} \theta_v(x_v) \span \label{eq:iLAP:a}\\
    \forall v\in\V:\;& x_v\in\L_v \\
    \forall \l\in\L\setminus\{\#\}:\;& \sum_{v\in\V} \iverson{x_v=\l}\leq1.\label{eq:iLAP:c}
\end{align}
\end{subequations}

\begin{remark}
As opposed to the LAP from Section~\ref{se:LAP}, we do not require~$|\V|=|\L|$ here. Next, the equality sign in~\eqref{eq:LAP:c} changed to the inequality sign in~\eqref{eq:iLAP:c}, i.e., not every label needs to be assigned to some vertex. Also, the constraint~\eqref{eq:iLAP:c} does not limit the number of vertices assigned to the dummy label. Thus, in contrast to the LAP, no question of feasibility arises because the assignment defined by~$x_v=\#$ for all~$v\in\V$ is always feasible.
\end{remark}

\begin{remark}
There exist many variants and extensions of the LAP~\cite{burkard2009assignment,volgenant1996linear,bijsterbosch2010solving,bertsekas1993reverse,chen2016conflict,murty1967symmetric,ramshaw2012minimum}. The closest problem to the ILAP is the \emph{rectangular LAP}~\cite{burkard2009assignment,bourgeois1971extension,bijsterbosch2010solving} (also called asymmetric~\cite{bertsekas1993reverse} or unbalanced LAP~\cite{ramshaw2012minimum}) where different sizes of the partitions are also assumed (with~$|\V|<|\L|$) but the dummy label is not considered because each vertex is required to be assigned to some non-dummy label. The \emph{coverage-sensitive many-to-many min-cost
bipartite matching} (CSM) problem from~\cite{lo2013evaluating} is more general than the ILAP as it allows specifying a (possibly infinite) cost for each label~$\L$ that depends on how many vertices are assigned this label.
\end{remark}

Analogously to the LAP, the ILAP can be also formulated as a linear program that is the left-hand problem of the primal-dual pair
\begin{subequations}\label{eq:iLAP_LP_relax}
\begin{align}
    && \min \; \sum_{\substack{v\in\V\\\l\in\L_v}} \theta_v(\l)\mu_v(\l) \span & \max \; \sum_{v\in\V} \alpha_v+\smashoperator{\sum_{\l\in\L\setminus\{\#\}}} \beta_\l \span\span\span \label{eq:iLAP_LP_relax:a}\\
    \forall v\in\V, \,\l\in\L_v:\;&& \mu_v(\l)&\geq 0 & \alpha_v+\iverson{\l\neq\#}\beta_\l&\leq\theta_v(\l)\label{eq:iLAP_LP_relax:b}\\
    \forall v\in\V:\;&& \sum_{\l\in\L_v} \mu_v(\l)&=1 & \alpha_v&\in\mathbb{R}\\
    \forall \l\in\L\setminus\{\#\}:\;&& \sum_{v\in\V_\l} \mu_v(\l)&\leq 1 & \beta_\l&\leq0.\label{eq:iLAP_LP_relax:d}
\end{align}
\end{subequations}
Again, we wrote the dual linear program on the right. Note that the dual variables~$\beta_\l$ are defined only for~$\l\in\L\setminus\{\#\}$, so the dual constraint~\eqref{eq:iLAP_LP_relax:b} does not contain~$\beta_{\l}$ if~$\l=\#$. $\V_\l$~is defined as in~\eqref{eq:LAP_LP_relax}.

It is not hard to show that the primal~\eqref{eq:iLAP_LP_relax} (on the left-hand side) is integral and that any its integral feasible solution corresponds to a feasible solution of~\eqref{eq:iLAP} with the same objective. This fact follows immediately, e.g., from~\cite[Theorem~5.21]{schrijver2003combinatorial} or~\cite[Theorem~4.2]{burkard2009assignment}.

\subsection{Reduction of the ILAP to the LAP}\label{se:reduction}

Next, we describe a reduction of the ILAP to the LAP that allows us to solve the ILAP instances and also obtain a relative-interior solution of the dual~\eqref{eq:iLAP_LP_relax} based on our previous results from Section~\ref{se:LAP}. A~similar reduction was informally mentioned in~\cite{thesisQAP} and also discussed in~\cite[Section~1.3]{ramshaw2012minimum} (which however transforms the rectangular LAP to the LAP instead of the ILAP to the LAP). The reduction provided in~\cite[Section~A1]{haller2022comparative} is also to some extent similar to the one considered here but it neither preserves the sparsity pattern nor splits the costs, so the resulting cost matrix is not symmetric.

\paragraph{Reduction}
For the purpose of showing our reduction formally, let~$\P=(\V,\L,\#,\{\L_v\}_{v\in\V},\theta)$ define an instance of the ILAP. We define the instance of LAP~$\P'=(\V',\L',\{\L_v'\}_{v\in \V'},\theta')$ by
\begin{subequations}\label{eq:transformation}
\begin{align}
    && \L'=\V'&=\V\cup\L\setminus\{\#\}\label{eq:transformation:a}\\
    \forall v\in \V':&&\L'_v &=\begin{cases}
        (\L_v\setminus\{\#\})\cup\{v\} & \text{ if } v\in \V\\
        \V_{v}\cup\{v\} & \text{ if } v \in \L
    \end{cases}\\
    \forall v\in \V',\,\l\in\L'_v:&&\theta'_v(\l) &= \begin{cases}
        \theta_v(\l)/2 & \text{ if } v\in \V,\, \l\in\L\\
        \theta_{\l}(v)/2 & \text{ if } v\in \L,\, \l\in\V\\
        \theta_{v}(\#) & \text{ if } v\in \V,\, \l\in\V\\
        0 & \text{ if } v\in \L,\, \l\in\L\\
    \end{cases}.\label{eq:transformation:c}
\end{align}
\end{subequations}
If a sparse encoding is used for the sets~$\L_v$, the size of the instance~$\P'$ scales linearly with the size of~$\P$. In this subsection, we analyze this reduction and assume that the instances~$\P$ and~$\P'$ are fixed for brevity.

\begin{remark}
Strictly speaking, we required in Section~\ref{se:LAP} the partitions~$\V'$ and~$\L'$ to be disjoint, which is not satisfied by~\eqref{eq:transformation:a}. However, there is no ambiguity since the sets of allowed labels~$\{\L'_v\}_{v\in\V'}$ are defined for elements of~$\V'$ and the costs~$\theta'$ are denoted asymetrically, i.e., for~$\theta'_v(\l)$, we have~$v\in\V'$ and~$\l\in\L'$ (and analogously for the primal variables~$\mu$ in~\eqref{eq:LAP_LP_relax}).  By distinguishing~$\V'$ from~$\L'$ in equations, we resolve any ambiguity.
\end{remark}

To provide the reader with an intuitive understanding of the reduction, let us now informally comment on it. For an assignment~$x\colon\V\rightarrow\L$ feasible for ILAP~$\P$, one can define an assignment~$x'\colon \V'\rightarrow\L'$ feasible for LAP~$\P'$ by
\begin{equation}\label{eq:perfect_matching_from_x}
    x'_{v'} = \begin{cases}
    x_{v'} &\text{ if } v'\in\V \text{ and } x_{v'}\neq\#\\
    v' &\text{ if } v'\in\V \text{ and } x_{v'}=\#\\
    u \text{ where~$u\in\V$ is such that~$x_u=v'$ }& \text{ if } v'\in\L \text{ and~$\exists u\in\V:x_u=v'$}\\
    v'& \text{ if } v'\in\L \text{ and~$\forall u\in\V:x_u\neq v'$}
    \end{cases}
\end{equation}
for all~$v'\in\V'$. In words, if~$v\in\V$ is assigned to~$x_v=\l\neq\#$, then~$x'_{v}=\l$ and~$x'_\l=v$. On the other hand, if~$v\in\V$ is assigned to the dummy label~$x_v=\#$, then it is assigned to itself by~$x'$, i.e., $x'_v=v$. If some label~$\l\in\L\setminus\{\#\}$ is not assigned to any vertex, i.e., $x_v\neq\l$ for all~$v\in\V$, then~$\l$ is also assigned to itself by~$x'$, i.e., $x'_\l=\l$.

By~\eqref{eq:transformation:c}, the cost for assigning a vertex~$v\in\V$ to itself in~$\P'$ is equal to assigning the vertex to the dummy label~$\#$ in~$\P$, i.e., $\theta_v(\#)$. The cost for assigning a label~$\l\in\L\setminus\{\#\}$ to itself in~$\P'$ is zero. Consequently, since the other costs are halved, the objective value of~$x'$ for LAP~$\P'$ is equal to the objective value of~$x$ for ILAP~$\P$. We prove later in Proposition~\ref{pr:same_opt_val} that the optimal values of~$\P$ and~$\P'$ are equal.

\begin{figure}[]
\centering
\begin{subfigure}[]{.35\textwidth}
    \centering
    \begin{tikzpicture}[baseline={(0,0)}]
        \node (pic) at (0,0) {\includegraphics[width=120pt]{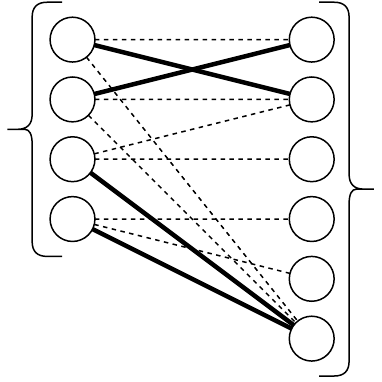}};
        \node[] at (-37.5pt,47.5pt) {\val{a}};
        \node[] at (-37.5pt,28.75pt) {\val{b}};
        \node[] at (-37.5pt,10pt) {\val{c}};
        \node[] at (-37.5pt,-8.75pt) {\val{d}};
        \node[] at (38pt,47.5pt) {\val{A}};
        \node[] at (38pt,28.75pt) {\val{B}};
        \node[] at (38pt,10pt) {\val{C}};
        \node[] at (38pt,-8.75pt) {\val{D}};
        \node[] at (38pt,-27.5pt) {\val{E}};
        \node[] at (38pt,-46.25pt) {\#};
        \node[] at (-65pt,19.25pt) {$\V$};
        \node[] at (64pt,-0pt) {$\L$}; 
    \end{tikzpicture}
    \caption{Bipartite graph with partitions~$\V$ and~$\L$ that defines an instance of the ILAP (up to the costs~$\theta$) and a feasible assignment~$x$ for this instance (determined by bold edges).}
    \label{fig:assignment1_ilap} 
\end{subfigure}\hspace{0.04\textwidth}
\begin{subfigure}[]{.275\textwidth}
    \centering
    \begin{tikzpicture}[baseline={(0,0)}]
        \node (pic) at (0,0) {\includegraphics[width=91.18pt]{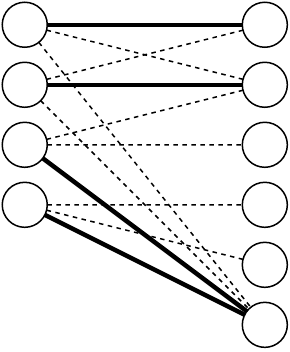}};
        \node[] at (-37.5pt,47.5pt) {\val{a}};
        \node[] at (-37.5pt,28.75pt) {\val{b}};
        \node[] at (-37.5pt,10pt) {\val{c}};
        \node[] at (-37.5pt,-8.75pt) {\val{d}};
        \node[] at (38pt,47.5pt) {\val{A}};
        \node[] at (38pt,28.75pt) {\val{B}};
        \node[] at (38pt,10pt) {\val{C}};
        \node[] at (38pt,-8.75pt) {\val{D}};
        \node[] at (38pt,-27.5pt) {\val{E}};
        \node[] at (38pt,-46.25pt) {\#};
    \end{tikzpicture}
    \caption{Assignment~$x^2$ feasible for the same ILAP instance as in Figure~\ref{fig:assignment1_ilap}.}
    \label{fig:assignment2_ilap} 
\end{subfigure}\hspace{0.04\textwidth}
\begin{subfigure}[]{.275\textwidth}
    \centering
    \begin{tikzpicture}[baseline={(0,0)}]
        \node (pic) at (0,0) {\includegraphics[width=91.18pt]{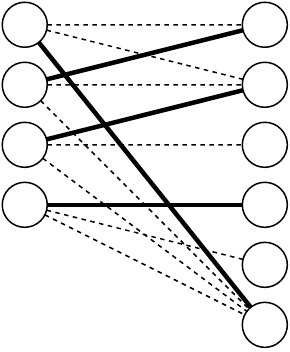}};
        \node[] at (-37.5pt,47.5pt) {\val{a}};
        \node[] at (-37.5pt,28.75pt) {\val{b}};
        \node[] at (-37.5pt,10pt) {\val{c}};
        \node[] at (-37.5pt,-8.75pt) {\val{d}};
        \node[] at (38pt,47.5pt) {\val{A}};
        \node[] at (38pt,28.75pt) {\val{B}};
        \node[] at (38pt,10pt) {\val{C}};
        \node[] at (38pt,-8.75pt) {\val{D}};
        \node[] at (38pt,-27.5pt) {\val{E}};
        \node[] at (38pt,-46.25pt) {\#};
    \end{tikzpicture}
    \caption{Assignment~$x^3$ feasible for the same ILAP instance as in Figure~\ref{fig:assignment1_ilap}.}
    \label{fig:assignment3_ilap} 
\end{subfigure}

\begin{subfigure}[]{.6\textwidth}
    \centering
    \begin{tikzpicture}[baseline={(0,0)}]
        \node (pic) at (0,0) {\includegraphics[width=218.74pt]{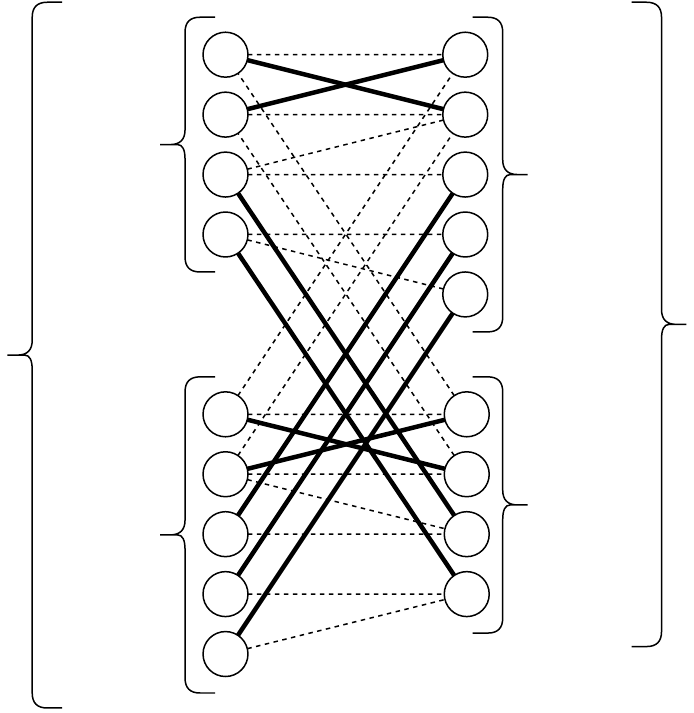}};
        \node[] at (-38.5pt,95pt) {\val{a}};
        \node[] at (-38.5pt,76.25pt) {\val{b}};
        \node[] at (-38.5pt,57.5pt) {\val{c}};
        \node[] at (-38.5pt,38.75pt) {\val{d}};

        \node[] at (37pt,95pt) {\val{A}};
        \node[] at (37pt,76.25pt) {\val{B}};
        \node[] at (37pt,57.5pt) {\val{C}};
        \node[] at (37pt,38.75pt) {\val{D}};
        \node[] at (37pt,20pt) {\val{E}};
        
        \node[] at (37.5pt,-74.75pt) {\val{d}};
        \node[] at (37.5pt,-56.0pt) {\val{c}};
        \node[] at (37.5pt,-37.25pt) {\val{b}};
        \node[] at (37.5pt,-18.5pt) {\val{a}};
        
        \node[] at (-38.5pt,-93.5pt) {\val{E}};
        \node[] at (-38.5pt,-74.75pt) {\val{D}};
        \node[] at (-38.5pt,-56.0pt) {\val{C}};
        \node[] at (-38.5pt,-37.25pt) {\val{B}};
        \node[] at (-38.5pt,-18.5pt) {\val{A}};
        
        \node[] at (-115pt,1pt) {$\V'$};
        \node[] at (115pt,10pt) {$\L'$};
        
        \node[] at (77pt,57pt) {$\L\setminus\{\#\}$};
        \node[] at (65pt,-47pt) {$\V$};
        \node[] at (-78pt,-56.25pt) {$\L\setminus\{\#\}$};
        \node[] at (-66pt,66.25pt) {$\V$};
    \end{tikzpicture}
    \caption{Instance of the LAP obtained by the reduction~\eqref{eq:transformation} from the ILAP instance determined by Figure~\ref{fig:assignment1_ilap}. The highlighted perfect matching~$x'$ is obtained by~\eqref{eq:perfect_matching_from_x} from the assignment in Figure~\ref{fig:assignment1_ilap}.}
    \label{fig:lap1} 
\end{subfigure}\hspace{0.04\textwidth}
\begin{subfigure}[]{.3\textwidth}
    \centering
    \begin{tikzpicture}[baseline={(0,0)}]
        \node (pic) at (-0.5pt,0) {\includegraphics[width=91.65pt]{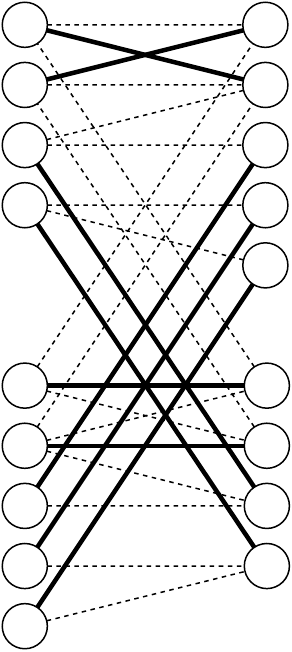}};
        \node[] at (-38.5pt,95pt) {\val{a}};
        \node[] at (-38.5pt,76.25pt) {\val{b}};
        \node[] at (-38.5pt,57.5pt) {\val{c}};
        \node[] at (-38.5pt,38.75pt) {\val{d}};

        \node[] at (37pt,95pt) {\val{A}};
        \node[] at (37pt,76.25pt) {\val{B}};
        \node[] at (37pt,57.5pt) {\val{C}};
        \node[] at (37pt,38.75pt) {\val{D}};
        \node[] at (37pt,20pt) {\val{E}};
        
        \node[] at (37.5pt,-74.75pt) {\val{d}};
        \node[] at (37.5pt,-56.0pt) {\val{c}};
        \node[] at (37.5pt,-37.25pt) {\val{b}};
        \node[] at (37.5pt,-18.5pt) {\val{a}};
        
        \node[] at (-38.5pt,-93.5pt) {\val{E}};
        \node[] at (-38.5pt,-74.75pt) {\val{D}};
        \node[] at (-38.5pt,-56.0pt) {\val{C}};
        \node[] at (-38.5pt,-37.25pt) {\val{B}};
        \node[] at (-38.5pt,-18.5pt) {\val{A}};
    \end{tikzpicture}
    \caption{A different assignment~$x''$ feasible for the LAP instance from Figure~\ref{fig:lap1}.}
    \label{fig:lap2} 
\end{subfigure}
\caption{Illustrations to Example~\ref{ex:transformation}.}
\label{fig:example_transformation}
\end{figure}
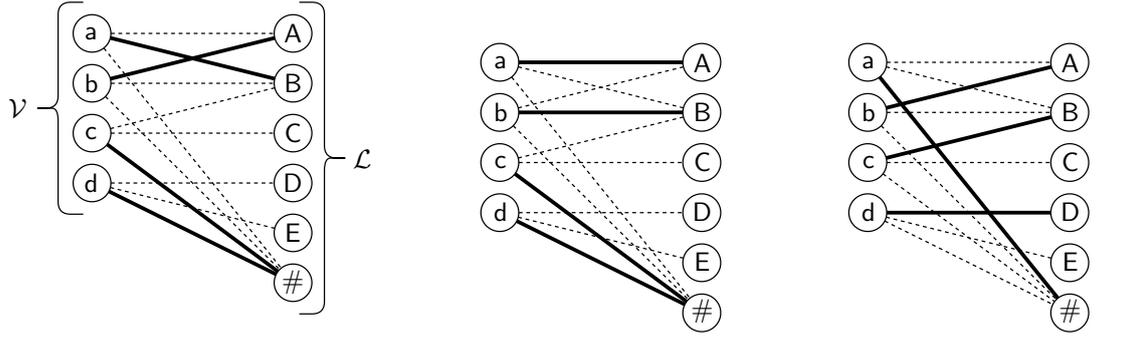
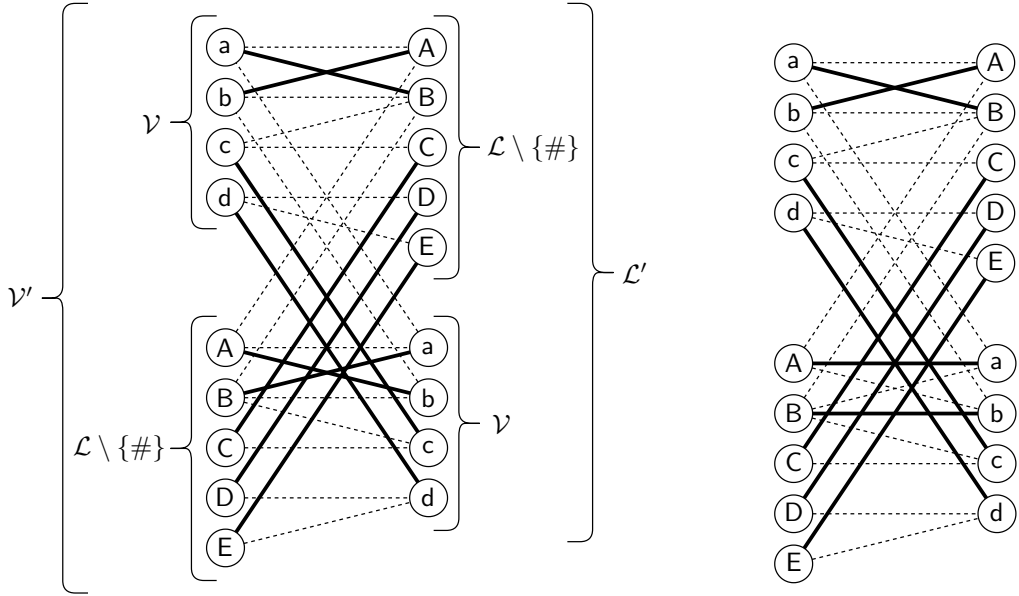

\begin{example}\label{ex:transformation}
To exemplify this reduction, let~$\V=\{\val{a},\val{b},\val{c},\val{d}\}$ and~$\L=\{\val{A},\val{B},\val{C},\val{D},\val{E},\#\}$. The allowed labels are given by~$\L_\val{a}=\L_\val{b}=\{\val{A},\val{B},\#\}$, $\L_\val{c}=\{\val{B},\val{C},\#\}$, and~$\L_\val{d}=\{\val{D},\val{E},\#\}$. This instance is depicted in Figure~\ref{fig:assignment1_ilap} where the allowed labels are connected by edges to the corresponding vertices.

The bold edges in Figure~\ref{fig:assignment1_ilap} define an assignment~$x\colon\V\rightarrow\L$ feasible for this instance, i.e., $x_\val{a}=\val{B}$, $x_\val{b}=\val{A}$, and~$x_\val{c}=x_\val{d}=\#$. Based on this instance~$\P$ of ILAP, one can define the LAP~$\P'$ by~\eqref{eq:transformation}, which is shown in Figure~\ref{fig:lap1}. Elements of~$\V'$ and~$\L'$ are on the left and right part of the bipartite graph, respectively. As in Figure~\ref{fig:assignment1_ilap}, the sets~$\{\L'_{v'}\}_{v'\in\V'}$ are depicted by drawing an edge~$\{v',\l'\}$ for each~$v'\in\V'$ and~$\l'\in \L'_{v'}$. We have, e.g., $\L'_\val{a}=\{\val{A},\val{B},\val{a}\}$ and~$\L'_\val{B}=\{\val{B},\val{a},\val{b},\val{c}\}$.  Figure~\ref{fig:lap1} also shows the assignment~$x'\colon \V'\rightarrow\L'$ defined by~\eqref{eq:perfect_matching_from_x} based on the aforementioned assignment~$x\colon\V\rightarrow\L$. This assignment corresponds to a perfect matching in the bipartite graph shown in Figure~\ref{fig:lap1}.

See that no vertex~$v\in\V$ is assigned to label~\val{C}, \val{D}, or~\val{E} in Figure~\ref{fig:assignment1_ilap}. Consequently, these labels are assigned to themselves in Figure~\ref{fig:lap1}. Similarly, vertices~\val{c} and~\val{d} are assigned to the dummy label~$\#$ in Figure~\ref{fig:assignment1_ilap} and thus, they are assigned to themselves in Figure~\ref{fig:lap1}. Informally speaking, for the other vertices and labels, the lower part of the diagram is obtained by mirroring the upper part horizontally.

Figure~\ref{fig:lap2} shows another assignment~$x''\colon\V'\rightarrow\L'$ feasible for this LAP instance. This assignment differs from~$x'$ only in the lower part of the LAP instance and the resulting perfect matching is thus not `symmetric'\,\footnote{More precisely, the bijection~$x'$ defined by Figure~\ref{fig:lap1} is an involution, i.e., it is self-inverse. The bijection~$x''$ defined by Figure~\ref{fig:lap2} is not an involution.}. From the lower part, one can extract the assignment~$x^2\colon\V\rightarrow\L$ that is shown in Figure~\ref{fig:assignment2_ilap} and is also feasible for the aforesaid ILAP instance.
\end{example}

Although one can map each assignment~$x$ feasible for the ILAP~$\P$ to an assignment~$x'$ feasible for the LAP~$\P'$ via~\eqref{eq:perfect_matching_from_x}, not every assignment~$x'$ feasible for the LAP~$\P'$ can be mapped to a unique assignment for the ILAP~$\P$. To be precise, one can generally extract two different ILAP assignments from each assignment feasible for the LAP. Recalling Example~\ref{ex:transformation}, this is the situation of the LAP assignment shown in Figure~\ref{fig:lap2} that can be decomposed into the two ILAP assignments in Figures~\ref{fig:assignment1_ilap} and~\ref{fig:assignment2_ilap}. Proposition~\ref{pr:2_assignments} analyzes the objective values of such decomposed assignments.

\begin{remark}\label{re:combining_2_assignments}
Seen from the other side, not every two assignments feasible for ILAP~$\P$ can be combined in this way to result in a feasible assignment for LAP~$\P'$. A sufficient condition for the possibility to combine assignments~$x,\,y\colon\V\rightarrow\L$ feasible for ILAP~$\P$ is that~$\{x_v\mid v\in\V\}=\{y_v\mid v\in\V\}$ and~$\{v\in\V\mid x_v=\# \}=\{v\in\V\mid y_v=\# \}$, i.e., if the sets of assigned labels are the same and the sets of vertices assigned to the dummy label are the same.

To illustrate this, let us continue in Example~\ref{ex:transformation}. The assignments~$x^2,\,x^3\colon\V\rightarrow\L$ defined by Figures~\ref{fig:assignment2_ilap} and~\ref{fig:assignment3_ilap}, respectively, do not satisfy these conditions because~$\{x^2_v\mid v \in \V\}=\{\val{A},\val{B},\#\}\neq\{x^3_v\mid v \in \V\}=\{\val{A},\val{B},\val{D},\#\}$ and even~$\{v\in\V\mid x_v^2=\# \}=\{\val{c},\val{d}\}\neq\{v\in\V\mid x^3_v=\# \}=\{\val{a}\}$. In this case, if the lower and upper part of the LAP~$\P'$ is assigned based on~$x^2$ and~$x^3$, respectively, such a matching cannot be completed to a feasible assignment for the LAP~$\P'$ (i.e., to a perfect matching in the bipartite graph shown in Figure~\ref{fig:lap1}).
\end{remark}


\begin{proposition}\label{pr:2_assignments}
Let~$x'$ be an assignment feasible for~$\P'$. Define assignments $x^1,x^2\colon\V\rightarrow\L$ by
\begin{subequations}\label{eq:def_2_assignments}
\begin{align}
    \forall v\in\V: x^1_v &= \begin{cases}
    x'_v & \text{ if } x'_v\in\L\setminus\{\#\}\\
    \# & \text{ otherwise (i.e., $x'_v=v$)}
    \end{cases}\\
    \forall v\in\V: x^2_v &= \begin{cases}
    x'^{-1}_v & \text{ if } x'^{-1}_v\in\L\setminus\{\#\}\\
    \# & \text{ otherwise (i.e., $x'^{-1}_v=v$) }
    \end{cases}
\end{align}
\end{subequations}
where~$x'^{-1}\colon\L'\rightarrow \V'$ is the inverse of~$x'$ (recall that~$x'$ is a bijection). Assignments~$x^1$ and~$x^2$ are feasible for the ILAP~$\P$. If~$\Theta'$ is the cost of~$x'$ for LAP~$\P'$ and~$\Theta_1$ and~$\Theta_2$ are the costs of~$x^1$ and~$x^2$ for ILAP~$\P$, respectively, then~$2\Theta'=\Theta_1+\Theta_2$.
\end{proposition}
\begin{proof}
To show feasibility of~$x^1$, note that~$x'_v$, $v\in\V$ are unique and that if~$x'_v\in\L$, then~$x'_v\in\L_v$. In detail, this is due to~$x'_v\in\L'_v=(\L_v\setminus\{\#\})\cup\{v\}$, so, if~$x'_v\notin\L$, then~$x'_v=v$. The proof for~$x^2$ is analogous since~$\V'_v=(\L_v\setminus\{\#\})\cup\{v\}$ and~$x'^{-1}_v\in\V'_v$ for all~$v\in\V$.

For the next statement, see that
\begin{subequations}
\begin{align}
    2\Theta' &= 2\sum_{v\in\V'}\theta'_v(x'_v) = 2\sum_{\substack{v\in\V\\x'_v\neq v}} \overbrace{\theta'_v(x'_v)}^{\theta_v(x'_v)/2}
    +\;2\sum_{\substack{v\in\V\\x'_v= v}} \overbrace{\theta'_v(x'_v)}^{\theta_v(\#)}
    +\;2\smashoperator{\sum_{\substack{\l\in\L\setminus\{\#\}\\x'_\l\neq \l}}} \overbrace{\theta'_\l(x'_\l)}^{\theta_{x'_\l}(\l)/2}
    +\;2\smashoperator{\sum_{\substack{\l\in\L\setminus\{\#\}\\x'_\l= \l}}} \overbrace{\theta'_\l(x'_\l)}^{0}\label{eq:theta_costs:a}\\
    &= \underbrace{\sum_{\substack{v\in\V\\x'_v\neq v}} \theta_v(x'_v)
    +\sum_{\substack{v\in\V\\x'_v= v}} \theta_v(\#)}_{\Theta_1}+\underbrace{\sum_{\substack{v\in\V\\x'^{-1}_v= v}} \theta_v(\#)
    +\sum_{\substack{v\in\V\\x'^{-1}_v\in\L\setminus\{\#\}}} \theta_{v}(x'^{-1}_v)}_{\Theta_2}\label{eq:theta_costs:b}
\end{align}
\end{subequations}
where equality~\eqref{eq:theta_costs:a} holds because the sets of indices in the individual sums form a partition of~$\V'$. The values in the upper brackets in~\eqref{eq:theta_costs:a} follow from the definition of~$\theta'$ in~\eqref{eq:transformation:c}. The equality~\eqref{eq:theta_costs:b} holds because~$x'_v=v\iff x'^{-1}_v=v$ and
\begin{equation}\label{eq:sets_equality}
 \{(x'_\l,\l)\mid \l\in\L\setminus\{\#\},\,x'_\l\neq\l\}=\{(v,x'^{-1}_v)\mid v\in\V, \,x'^{-1}_v\in\L\setminus\{\#\}\}   
\end{equation}
which is due to~$\L'_\l=\V_\l\cup\{\l\}$ for all~$\l\in\L\setminus\{\#\}$, so the term~$x'_\l\neq\l$ in the first bracket in~\eqref{eq:sets_equality} is equivalent to~$x'_\l\in\V$. By~\eqref{eq:def_2_assignments}, the parts enclosed in the brackets in~\eqref{eq:theta_costs:b} correspond to~$\Theta_1$ and~$\Theta_2$.
\end{proof}

\paragraph{Properties of Optimal Solutions}

\begin{proposition}\label{pr:same_opt_val}
The optimal values of ILAP~$\P$ and LAP~$\P'$ are equal.
\end{proposition}
\begin{proof}
We prove this claim by showing that, for any assignment feasible for~$\P$, we can construct an assignment feasible for~$\P'$ that has the same or better objective and vice versa.

As discussed above, for any assignment~$x$ feasible for the ILAP~$\P$, we can construct an assignment~$x'$ via~\eqref{eq:transformation} that is feasible for LAP~$\P'$ and has the same objective.

For the other direction, let~$x'$ be a feasible assignment for the LAP~$\P'$ with cost~$\Theta'$. By Proposition~\ref{pr:2_assignments}, one can construct the assignments~$x^1$ and~$x^2$ whose objective values are~$\Theta_1$ and~$\Theta_2$, respectively. Without loss of generality, let~$\Theta_1\leq\Theta_2$. Then, $\Theta'=(\Theta_1+\Theta_2)/2\geq \Theta_1$ and assignment~$x^1$ is feasible for~$\P$ and has the same or better objective than~$x'$ for~$\P'$.
\end{proof}

\begin{remark}
Let~$x^1$ and~$x^2$ be obtained from~$x'$ as in Proposition~\ref{pr:2_assignments}. It follows from Proposition~\ref{pr:2_assignments} (combined with Proposition~\ref{pr:same_opt_val}) that the assignment~$x'$ is optimal for LAP~$\P'$ if and only if both~$x^1$ and~$x^2$ are optimal for ILAP~$\P$. Consequently, if ILAP~$\P$ has~$n$ optimal solutions, then LAP~$\P'$ has at most~$n^2$ optimal solutions. Note that LAP~$\P'$ need not have exactly~$n^2$ optimal solutions as not every pair of optimal solutions of ILAP~$\P$ can be combined together (recall Remark~\ref{re:combining_2_assignments}).

To compare with other reductions connected to the LAP, e.g., the well-known reduction from the rectangular LAP to the LAP in~\cite[Section~5.4.4]{burkard2009assignment} (also mentioned in~\cite{bijsterbosch2010solving}) increases the overall number of optimal solutions $(|\L|-|\V|)!$ times where~$|\V|$ and~$|\L|$ are the partition sizes of the rectangular LAP (with~$|\V|<|\L|$)~\cite{bijsterbosch2010solving}. The reduction of the IQAP to the QAP in~\cite[Section~A1]{haller2022comparative} (and the analogous reduction of the ILAP to the LAP) increases the number of optimal solutions at least ${\max\{|\L|,|\V|\}!}/{||\L|-|\V||!}$ times.
\end{remark}

Now, we will discuss how to transform a solution feasible for the LP formulation of LAP~$\P'$ to a solution feasible for the LP formulation of ILAP~$\P$. The mappings for both primal and dual solutions are introduced in Theorem~\ref{th:maps_ri_to_ri} where we show that they preserve not only optimality but also membership in the relative interior of optimizers. For the proof of this theorem, we require an auxiliary lemma.

\begin{lemma}\label{le:rel_int_symmetric}
Let~$\mu'$ be from the relative interior of optimizers of the primal~\eqref{eq:LAP_LP_relax} for LAP~$\P'$ (i.e., from the relative interior of optimizers of the primal~\eqref{eq:transformation_LP}, stated later). Then, for all~$v\in\V$ and~$\l\in\L_v\setminus\{\#\}$, $\mu'_v(\l)>0\iff\mu'_\l(v)>0$.
\end{lemma}
\begin{proof}
Let~$\mu^*$ be defined by
\begin{subequations}
\begin{alignat}
    \forall v\in\V:\;&&\mu^*_v(v) & = \mu'_v(v)\\
    \forall \l\in\L\setminus\{\#\}:\;&&\mu^*_\l(\l) & = \mu'_\l(\l)\\
    \forall v\in\V,\,\l\in\L_v\setminus\{\#\}:\;&& \mu^*_v(\l) &=\mu'_\l(v)\\
    \forall v\in\V,\,\l\in\L_v\setminus\{\#\}:\;&& \mu^*_\l(v) &=\mu'_v(\l).
\end{alignat}
\end{subequations}
By symmetry in the definition of~$\mu^*$, it is easy to see that $\mu^*$ is also feasible for the primal~\eqref{eq:transformation_LP}. Moreover, the objective values coincide for~$\mu^*$ and~$\mu'$, so~$\mu^*$ is also optimal.

We prove the claim by contradiction.\footnote{Lemma~\ref{le:rel_int_symmetric} can be also proved using Proposition~\ref{pr:rel_int_primal}: by symmetry of the problem~$\P'$, there is an optimal assignment with~$x_v=\l$ if and only if there is an optimal assignment with~$x_\l=v$.} Let~$(\alpha',\beta')$ be from the relative interior of optimizers of the dual~\eqref{eq:transformation_LP}. Without loss of generality, suppose that there is~$v\in\V$ and~$\l\in\L_v\setminus\{\#\}$ such that~$\mu'_v(\l)>0$ and~$\mu'_\l(v)=0$. By strict complementary slackness, $\alpha'_\l+\beta'_v<\theta_v(\l)/2$ (see~\eqref{eq:transformation_LP:e}). By definition of~$\mu^*$, we have~$\mu^*_\l(v)=\mu'_v(\l)>0$ and~$\mu^*$ does not satisfy complementary slackness with~$(\alpha',\beta')$. This is contradictory with optimality of~$\mu^*$.
\end{proof}

\begin{theorem}\label{th:maps_ri_to_ri}
Let~$\mu'$ and~$(\alpha',\beta')$ be feasible for the primal and dual~\eqref{eq:LAP_LP_relax}, respectively, for the LAP instance~$\P'$. Define~$\mu$ by 
\begin{equation}\label{eq:mu_def}
    \forall v\in\V,\,\l\in\L_v: \mu_v(\l)=\begin{cases}
    (\mu'_v(\l)+\mu'_\l(v))/2 &\text{ if } \l\neq\#\\
    \mu'_v(v) &\text{ if } \l=\#
    \end{cases}
\end{equation}
and~$(\alpha,\beta)$ by
\begin{subequations}\label{eq:alpha_beta}
\begin{align}
    \forall v\in\V: \, \alpha_v &= \alpha'_v+\beta'_v\\
    \forall \l\in\L\setminus\{\#\}: \, \mathrlap{\beta_{\l}}\phantom{\alpha_v} &= \alpha'_\l+\beta'_\l.
\end{align}
\end{subequations}
Then it holds that
\begin{enumerate}
    \item[(a)] $\mu$~is feasible for the primal~\eqref{eq:iLAP_LP_relax},
    \item[(b)] $(\alpha,\beta)$~is feasible for the dual~\eqref{eq:iLAP_LP_relax},
    \item[(c)] if~$\mu'$~is optimal for primal~\eqref{eq:LAP_LP_relax}, then $\mu$~is optimal for primal~\eqref{eq:iLAP_LP_relax},
    \item[(d)] if~$(\alpha',\beta')$~is optimal for dual~\eqref{eq:LAP_LP_relax}, then $(\alpha,\beta)$~is optimal for dual~\eqref{eq:iLAP_LP_relax},
    \item[(e)] if~$\mu'$~is in the relative interior of optimizers of primal~\eqref{eq:LAP_LP_relax}, then $\mu$~is in the relative interior of optimizers of primal~\eqref{eq:iLAP_LP_relax},
    \item[(f)] if~$(\alpha',\beta')$~is in the relative interior of optimizers of dual~\eqref{eq:LAP_LP_relax}, then $(\alpha,\beta)$~is in the relative interior of optimizers of dual~\eqref{eq:iLAP_LP_relax}.
\end{enumerate}
\end{theorem}
\begin{proof}
For clarity of this proof, we write the LP formulation~\eqref{eq:LAP_LP_relax} of the LAP~$\P'$. Together with the corresponding dual, this reads
\begin{subequations}\label{eq:transformation_LP}
\begin{align}
    && \min \; p(\mu')\span & \max\; d(\alpha',\beta') \span\span\span\\
    \forall v \in \V, \,\l\in\L_v\setminus\{\#\}:&& \mu'_v(\l)&\geq0 & \alpha'_v+\beta'_\l&\leq \theta_v(\l)/2 \label{eq:transformation_LP:b}\\
    \forall v \in \V:&& \mu'_v(v)&\geq0 & \alpha'_v+\beta'_v&\leq \theta_v(\#)\label{eq:transformation_LP:c}\\
    \forall \l \in \L\setminus\{\#\}:&& \mu'_\l(\l)&\geq0 & \alpha'_\l+\beta'_\l&\leq 0\label{eq:transformation_LP:d}\\
    \forall \l\in \L\setminus\{\#\},\,v\in\V_\l:&& \mu'_\l(v)&\geq0&\alpha'_\l+\beta'_v&\leq \theta_v(\l)/2\label{eq:transformation_LP:e}\\
    \forall v\in\V:&& \mu'_v(v)+\smashoperator{\sum_{\l\in\L_v\setminus\{\#\}}}\mu'_v(\l)&=1& \alpha'_v&\in\mathbb{R}\label{eq:transformation_LP:f}\\
    \forall \l\in\L\setminus\{\#\}:&& \mu'_\l(\l)+\sum_{v\in \V_\l} \mu'_\l(v) &=1& \alpha'_\l&\in\mathbb{R}\label{eq:transformation_LP:g}\\
    \forall v\in\V:&& \mu'_v(v)+\smashoperator{\sum_{\l\in\L_v\setminus\{\#\}}}\mu'_\l(v)&=1& \beta'_v&\in\mathbb{R}\label{eq:transformation_LP:h}\\
    \forall \l\in\L\setminus\{\#\}:&& \mu'_\l(\l)+\sum_{v\in\V_\l} \mu'_v(\l)&=1& \beta'_\l&\in\mathbb{R}\label{eq:transformation_LP:i}
\end{align}
\end{subequations}
where the primal and dual objectives are defined by
\begin{subequations}\label{eq:objectives}
\begin{align}
    p(\mu')&= {\sum_{\substack{v\in\V\\\l\in\L\setminus\{\#\}}}} \theta_v(\l)\big(\mu'_v(\l)+\mu'_\l(v) \big)/2 + \sum_{v\in\V} \theta_v(\#)\mu'_v(v)\\
    d(\alpha',\beta')&= \sum_{v\in\V}(\alpha'_v+\beta'_v) + {\sum_{\l\in\L\setminus\{\#\}}} (\alpha'_\l+\beta'_\l).
\end{align}
\end{subequations}
Note that lines~\eqref{eq:transformation_LP:b}-\eqref{eq:transformation_LP:e}, \eqref{eq:transformation_LP:f}-\eqref{eq:transformation_LP:g}, and~\eqref{eq:transformation_LP:h}-\eqref{eq:transformation_LP:i} correspond to lines~\eqref{eq:LAP_LP_relax:b}, \eqref{eq:LAP_LP_relax:c}, and~\eqref{eq:LAP_LP_relax:d}, respectively.

We begin by proving~(a). For any~$v\in\V$, it holds that
\begin{equation}
    \sum_{\l\in\L_v}\mu_v(\l) =
    \Big(\mu'_v(v)+{\sum_{\l\in\L_v\setminus\{\#\}}}\mu'_v(\l)\Big)/2 +
    \Big(\mu'_v(v)+{\sum_{\l\in\L_v\setminus\{\#\}}}\mu'_\l(v)\Big)/2 = 1
\end{equation}
where we used~\eqref{eq:mu_def} and primal constraints~\eqref{eq:transformation_LP:f} and~\eqref{eq:transformation_LP:h}. Analogously, for any~$\l\in\L\setminus\{\#\}$,
\begin{equation}
    \sum_{v\in\V_\l} \mu_v(\l) = \sum_{v\in\V_\l}\mu'_v(\l)/2+\sum_{v\in\V_\l}\mu'_\l(v)/2 \leq 1
\end{equation}
which again follows from~\eqref{eq:mu_def} and primal constraints~\eqref{eq:transformation_LP:g} and~\eqref{eq:transformation_LP:i}. Non-negativity of~$\mu$ follows trivially from non-negativity of~$\mu'$.

Next, we proceed with~(b). To show that~$(\alpha,\beta)$ satisfies the dual constraint~\eqref{eq:iLAP_LP_relax:b}, we consider two cases, depending on whether~$\l=\#$. First, if~$v\in\V$ and~$\l\in\L_v\setminus\{\#\}$, then~$\alpha_v+\beta_\l=\alpha'_v+\beta'_v+\alpha'_\l+\beta'_\l\leq\theta_v(\l)$ where the equality is given by~\eqref{eq:alpha_beta} and the inequality follows from the two corresponding inequalities in~\eqref{eq:transformation_LP:b} and~\eqref{eq:transformation_LP:e}. Second, if~$v\in\V$ and~$\l=\#$, then~$\alpha_v=\alpha'_v+\beta'_v\leq\theta_v(\#)$ where the equality again follows from~\eqref{eq:alpha_beta} and the inequality from~\eqref{eq:transformation_LP:c}. Finally, for the other dual constraint~\eqref{eq:iLAP_LP_relax:d}, we have~$\beta_\l=\alpha'_\l+\beta'_\l\leq0$ by~\eqref{eq:alpha_beta} and~\eqref{eq:transformation_LP:d}.

To show~(c) and~(d), recall that the optimal values of ILAP~$\P$ and LAP~$\P'$ and their LP formulations coincide by Proposition~\ref{pr:same_opt_val} and integrality of the primal linear programs. By~(a) and~(b), \eqref{eq:mu_def}~and~\eqref{eq:alpha_beta} map a feasible solution of the primal and dual~\eqref{eq:transformation_LP} to a feasible solution of the primal and dual~\eqref{eq:iLAP_LP_relax}, respectively. The claim follows because these mappings preserve the objective, which can be verified by plugging the definitions~\eqref{eq:mu_def} and~\eqref{eq:alpha_beta} to~\eqref{eq:objectives} and comparing to~\eqref{eq:iLAP_LP_relax:a}.

To prove~(e) and~(f), let~$\mu'$ and~$(\alpha',\beta')$ be from the relative interior of optimizers of the primal and dual~\eqref{eq:transformation_LP}, respectively. We will prove that~$\mu$ defined by~\eqref{eq:mu_def} and~$(\alpha,\beta)$ defined by~\eqref{eq:alpha_beta} are in the relative interior of optimizers of the primal and dual~\eqref{eq:iLAP_LP_relax}, respectively, by showing that they satisfy strict complementary slackness.

We start by proving it for the constraints~\eqref{eq:iLAP_LP_relax:d}. See that, for all~$\l\in\L\setminus\{\#\}$,
\begin{equation}\label{eq:equiv_1}
        \beta_\l = 0 \iff \alpha'_\l+\beta'_\l=0 \iff \mu'_\l(\l)>0 \iff \sum_{v\in\V_\l} \mu'_v(\l) +\sum_{v\in \V_\l} \mu'_\l(v)  <2 \iff  \sum_{v\in\V_\l} \mu_v(\l)< 1
\end{equation}
where the first equivalence is given by definition of~$\beta$ in~\eqref{eq:alpha_beta}, the second equivalence follows from strict complementarity of~$\mu'$ and~$(\alpha',\beta')$ (see~\eqref{eq:transformation_LP:d}), and the third equivalence follows from constraints~\eqref{eq:transformation_LP:g} and~\eqref{eq:transformation_LP:i}. Definition of~$\mu$ in~\eqref{eq:mu_def} yields the last equivalence in~\eqref{eq:equiv_1}.

Next, we proceed with the constraints~\eqref{eq:iLAP_LP_relax:b}. For all~$v\in\V$ and~$\l\in\L_v\setminus\{\#\}$, the following equivalences hold:
\begin{equation}\label{eq:equiv_2}
 \alpha'_\l+\beta'_v=\theta_v(\l)/2 \iff \mu'_\l(v)> 0 \iff \mu'_v(\l) > 0 \iff \alpha'_v+\beta'_\l=\theta_v(\l)/2.
\end{equation}
The first and last equivalences in~\eqref{eq:equiv_2} are given by strict complementarity of~$\mu'$ and~$(\alpha',\beta')$ (see~\eqref{eq:transformation_LP:e} and~\eqref{eq:transformation_LP:b}) and the middle equivalence follows from Lemma~\ref{le:rel_int_symmetric}. These equivalences imply~$\mu_v(\l)>0 \iff \alpha_v+\beta_\l=\theta_v(\l)$ by definition of~$\mu$ and~$(\alpha,\beta)$ in~\eqref{eq:mu_def} and~\eqref{eq:alpha_beta}, respectively. For~$v\in\V$ and~$\l=\#$, we have that
\begin{equation}
    \mu_v(\#) > 0 \iff \mu'_v(v)>0 \iff \alpha'_v+\beta'_v=\theta_v(\#) \iff \alpha_v=\theta_v(\#)
\end{equation}
where the first and last equivalences hold by definition of~$\mu$ and~$(\alpha,\beta)$ and the middle equivalence holds by strict complementarity of~$\mu'$ and~$(\alpha',\beta')$ (see~\eqref{eq:transformation_LP:c}).
\end{proof}

\section{The Incomplete QAP and Subproblems of Its LP Relaxation}\label{se:iQAP}

We define the \emph{incomplete quadratic assignment problem} (IQAP)~\cite{hutschenreiter2021fusion,haller2022comparative} as follows. Let~$(\V,\E)$ be an undirected loopless graph where~$\V$ is a finite set of vertices and~$\E\subseteq \{\{u,v\}\mid u,v\in\V,\,u\neq v\}$ is a set of edges. For clarity of notation, we abbreviate~$\{u,v\}$ to~$uv$ in this section. Next, let $\L$~be a finite set of labels such that~$\#\in\L$ and, for each~$v\in\V,$ let $\L_v\subseteq \L$ with~$\#\in\L_v$ be the set of allowed labels for vertex~$v$. Finally, let~$\theta_v\colon \L_v\rightarrow \mathbb{R}$ and $\theta_{uv}\colon \L_u\times \L_v\rightarrow \mathbb{R}$ be the cost functions for each~$v\in \V$ and~$uv\in\E$, respectively (adopting that~$\theta_{uv}(\l,\k)=\theta_{vu}(\k,\l)$).

In this setting, the IQAP reads
\begin{subequations}\label{eq:iQAP}
\begin{align}
    \min \sum_{v\in\V} \theta_v(x_v)+\sum_{uv\in\E} \theta_{uv}(x_u,x_v) \span\label{eq:iQAP:a} \\
    \forall v\in\V:\;& x_v\in\L_v \\
    \forall \l\in\L\setminus\{\#\}:\;& \sum_{v\in\V} \iverson{x_v=\l}\leq 1\label{eq:iQAP:c}
\end{align}
\end{subequations}
and differs from the ILAP~\eqref{eq:iLAP} only by the quadratic term in the objective~\eqref{eq:iQAP:a}.

\begin{remark}
To compare, the (complete) QAP does not generally require a dummy label (i.e., we remove~$\#$ from~$\L$ and from each~$\L_v,\,v\in\V$) and seeks an assignment~$x\colon\V\rightarrow\L$ minimizing~\eqref{eq:iQAP:a} such that~$x_u\neq x_v$ for each distinct $u,\,v\in\V$. In this case, we typically also have~$|\V|=|\L|$ so that~$x$ is a bijection between~$\V$ and~$\L$~\cite{zhang2016pairwise,burkard2009assignment}.
\end{remark}

An LP relaxation of the IQAP is the left-hand problem of the primal-dual pair
\begin{subequations}\label{eq:iQAP_LP_relax}
\begin{align}
    \hspace{1.5cm}\min \; \sum_{\substack{v\in\V\\\l\in\L_v}} \theta_v(\l)\mu_v(\l)+\smashoperator{\sum_{\substack{uv\in\E\\(\k,\l)\in\L_u\times\L_v}}} \theta_{uv}(\k,\l) \mu_{uv}(\k,\l) \span\span\span & \max \; \sum_{v\in\V} \alpha_v+\smashoperator{\sum_{\l\in\L\setminus\{\#\}}} \beta_\l \span \label{eq:iQAP_LP_relax:a}\\
    \forall v\in\V, \,\l\in\L_v:\;&& \mu_v(\l)&\geq 0 & \alpha_v+\iverson{\l\neq\#}\beta_\l&\leq\theta_v^\phi(\l)\label{eq:iQAP_LP_relax:b}\\
    \forall uv\in\E,\,(\k,\l)\in\L_u\times\L_v:\;&& \mu_{uv}(\k,\l)&\geq 0 & 0&\leq \theta^\phi_{uv}(\k,\l)\\
    \forall v\in\V:\;&& \sum_{\l\in\L_v} \mu_v(\l)&=1 & \alpha_v&\in\mathbb{R}\\
    \forall v\in\V,\, u\in \N_v,\,\l\in\L_v:\;&& \sum_{\k\in\L_u} \mu_{uv}(\k,\l)&=\mu_v(\l) & \phi_{v\rightarrow u}(\l)&\in\mathbb{R}\\
    \forall \l\in\L\setminus\{\#\}:\;&& \sum_{v\in\V_\l} \mu_v(\l)&\leq 1 & \beta_\l&\leq 0
\end{align}
\end{subequations}
where~$\N_v=\{u\in \V\mid uv\in\E\}$ is the set of neighbors of vertex~$v$ in the graph~$(\V,\E)$ and
\begin{subequations}
\begin{align}
    \theta^\phi_v(\l)&=\theta_v(\l)+\sum_{u\in\N_v}\phi_{v\rightarrow u}(\l)\\
    \theta^\phi_{uv}(\k,\l)&=\theta_{uv}(\k,\l)-\phi_{v\rightarrow u}(\l)-\phi_{u\rightarrow v}(\k)
\end{align}
\end{subequations}
are reparametrized costs~\cite{savchynskyy2019discrete,hutschenreiter2021fusion}. In contrast to the problems considered in the previous sections, the primal~\eqref{eq:iQAP_LP_relax} is not integral and it is therefore not an LP formulation of the IQAP. This is of course not surprising as the IQAP problem is NP-hard.

The LP relaxation~\eqref{eq:iQAP_LP_relax} is implicit in~\cite[Section~5.2]{hutschenreiter2021fusion} and can be seen as an adaptation of the LP relaxation of the QAP from~\cite{zhang2016pairwise} to the case of the IQAP (also see~\cite[Section~3.3.1]{thesisQAP}). From a different point of view, the primal~\eqref{eq:iQAP_LP_relax} is an asymmetric analog of the classical Adams-Johnson linearization~\cite{adams1994improved}. The latter contains symmetric constraints with vertices and labels being swapped.

As already noted earlier, the beneficial feature of the relaxation~\eqref{eq:iQAP_LP_relax} is that block-variable subproblems of the dual correspond to dual LP formulation of the ILAP and dual LP relaxation of the WCSP. We thoroughly show this in the following two subsections.

\subsection{The ILAP Subproblem}

Let us fix the~$\phi$ variables in the dual~\eqref{eq:iQAP_LP_relax}. The dual restricted to the variables~$(\alpha,\beta)$ together with the corresponding primal\footnote{Since the primal-dual pair~\eqref{eq:iLAP_LP_relax_subproblem} is defined by restricting the set of dual variables in~\eqref{eq:iQAP_LP_relax}, the primal~\eqref{eq:iLAP_LP_relax_subproblem} contains only a subset of the primal constraints from~\eqref{eq:iQAP_LP_relax}. This also applies to the LP relaxation of the WCSP~\eqref{eq:MRF_LP_relax} shown in Section~\ref{se:WCSP_subproblem}.} reads
\begin{subequations}\label{eq:iLAP_LP_relax_subproblem}
\begin{align}
    && \min \; \sum_{\substack{v\in\V\\\l\in\L_v}} \theta_v^\phi(\l)\mu_v(\l) \span & \max \; \sum_{v\in\V} \alpha_v+\smashoperator{\sum_{\l\in\L\setminus\{\#\}}} \beta_\l \span\span\span \\
    \forall v\in\V, \,\l\in\L_v:\;&& \mu_v(\l)&\geq 0 & \alpha_v+\iverson{\l\neq\#}\beta_\l&\leq\theta_v^\phi(\l)\label{eq:iLAP_LP_relax_subproblem:b}\\
    \forall v\in\V:\;&& \sum_{\l\in\L_v} \mu_v(\l)&=1 & \alpha_v&\in\mathbb{R}\\
    \forall \l\in\L\setminus\{\#\}:\;&& \sum_{v\in\V_\l} \mu_v(\l)&\leq 1 & \beta_\l&\leq 0.
\end{align}
\end{subequations}
It is easy to see that this is the LP formulation~\eqref{eq:iLAP_LP_relax} of the ILAP except that the costs are~$\theta^\phi$ instead of~$\theta$.

For the purpose of applying BCA to the dual~\eqref{eq:iLAP_LP_relax_subproblem}, it is convenient to notice that, at optimum of the dual, we always have
\begin{equation}\label{eq:optimal_alpha}
    \forall v\in\V: \alpha_v=\min_{\l\in\L_v}\big( \theta_v^\phi(\l)-\iverson{\l\neq\#}\beta_\l\big).
\end{equation}
Plugging~\eqref{eq:optimal_alpha} into the dual~\eqref{eq:iLAP_LP_relax_subproblem} results in the optimization problem
\begin{subequations}\label{eq:iLAP_cpaf}
\begin{align}
    \max \; \sum_{v\in\V} &\min_{\l\in\L_v}\big( \theta_v^\phi(\l)-\iverson{\l\neq\#}\beta_\l\big)+\smashoperator{\sum_{\l\in\L\setminus\{\#\}}} \beta_\l  \label{eq:iLAP_cpaf:a}\\
    &\forall \l\in\L\setminus\{\#\}:\beta_\l\leq 0,
\end{align}
\end{subequations}
which can be interpreted as the maximization of a concave piecewise-affine function over non-positive variables.

\paragraph{Coordinate Ascent} To optimize~\eqref{eq:iLAP_cpaf} coordinate-wise, we derive coordinate-ascent updates for the individual~$\beta$ variables. For this, let~$\l\in\L\setminus\{\#\}$. The objective~\eqref{eq:iLAP_cpaf:a} restricted to the variable~$\beta_\l$ reads (up to a constant)
\begin{equation}\label{eq:obj_restricted}
    \sum_{v\in\V_\l} \min\left\{\theta_v^\phi(\l)-\beta_\l,\Phi_v(\l)\right\} + \beta_\l
\end{equation}
where
\begin{equation}
    \Phi_v(\l) = \min_{\l'\in\L_v\setminus \{\l\}}\big(\theta_v^\phi(\l')-\iverson{\l'\neq\#}\beta_{\l'}\big)
\end{equation}
are constants w.r.t.~$\beta_\l$. Let~$b_1$ and~$b_2$ be the smallest and second smallest value among~$\theta_v^\phi(\l)-\Phi_v(\l)$ for $v\in\V_\l$, respectively.\footnote{For clarity, if the minimal value is attained for multiple~$v\in\V_\l$, then~$b_1=b_2$. If~$|\V_\l|=1$ or~$\V_\l=\emptyset$, then~$b_2=0$ or~$b_1=b_2=0$, respectively. Note that, if~$|\V_\l|\geq 2$, then the objective~\eqref{eq:obj_restricted} is increasing for~$\beta_\l\leq b_1$, constant for~$b_1\leq\beta_\l\leq b_2$, and decreasing for~$\beta_\l\geq b_2$.} Based on~\cite[Lemma~28 in supplement]{werner2020relative}, the set of maximizers of~\eqref{eq:obj_restricted} subject to~$\beta_\l\leq 0$ is the interval~$[\min\{b_1,0\},\min\{b_2,0\}]$. To satisfy the relative-interior rule, we can choose, e.g., the midpoint of this interval, i.e., set
\begin{equation}\label{eq:coordinate_ascent}
    \beta_\l:=(\min\{b_1,0\}+\min\{b_2,0\})/2.
\end{equation}

\subsection{The WCSP Subproblem}\label{se:WCSP_subproblem}

Analogously, let us fix variables~$\beta$ in the dual LP relaxation~\eqref{eq:iQAP_LP_relax}. The dual restricted to variables~$(\alpha,\phi)$ together with the corresponding primal reads
\begin{subequations}\label{eq:MRF_LP_relax}
\begin{align}
   \hspace{1.5cm} \span\span\span \min \; \sum_{\substack{v\in\V\\\l\in\L_v}} \theta'_v(\l)\mu_v(\l)+\smashoperator{\sum_{\substack{uv\in\E\\(\k,\l)\in\L_u\times\L_v}}} \theta_{uv}(\k,\l) \mu_{uv}(\k,\l)  & \span \max \; \sum_{v\in\V} \alpha_v \\
    \forall v\in\V, \,\l\in\L_v:\;&& \mu_v(\l)&\geq 0 & \alpha_v&\leq\theta_v'^\phi(\l)\label{eq:MRF_LP_relax:b}\\
    \forall uv\in\E,\,(\k,\l)\in\L_u\times\L_v:\;&& \mu_{uv}(\k,\l)&\geq 0 & 0&\leq \theta^\phi_{uv}(\k,\l)\\
    \forall v\in\V:\;&& \sum_{\l\in\L_v} \mu_v(\l)&=1 & \alpha_v&\in\mathbb{R}\label{eq:MRF_LP_relax:d}\\
    \forall v\in\V,\, u\in \N_v,\,\l\in\L_v:\;&& \sum_{\k\in\L_u} \mu_{uv}(\k,\l)&=\mu_v(\l) & \phi_{v\rightarrow u}(\l)&\in\mathbb{R}\label{eq:MRF_LP_relax:e}
\end{align}
\end{subequations}
where~$\theta'_v(\l)=\theta_v(\l)-\iverson{\l\neq\#}\beta_\l$ are constants for each~$v\in\V,\,\l\in\L_v$. This subproblem is equivalent to the LP relaxation of the (pairwise) WCSP which is also known as the MAP inference problem in graphical models~\cite{zhang2016pairwise,hutschenreiter2021fusion,savchynskyy2019discrete,werner2007linear}. The LP relaxation~\eqref{eq:MRF_LP_relax} was proposed independently multiple times (e.g., in~\cite{schlesinger1976sintaksicheskiy} or~\cite{cooper2007optimal}) and is often referred to as the basic LP relaxation of the WCSP~\cite{thapper2012power} or local polytope relaxation~\cite{savchynskyy2019discrete}. Note that this linear program was shown to be as hard to solve as any linear program~\cite{prusa2013universality}, thus solving~\eqref{eq:iQAP_LP_relax} is also at least as hard.

\begin{remark}
There exist several equivalent formulations of the LP relaxation~\eqref{eq:MRF_LP_relax}~\cite{werner2007linear}. As with~\eqref{eq:iLAP_LP_relax_subproblem}, for any dual-optimal solution~$(\alpha,\phi)$ of~\eqref{eq:MRF_LP_relax}, \eqref{eq:optimal_alpha} holds, which results in a different form of the dual, namely
\begin{subequations}\label{eq:dual_MRF_LP_relax}
\begin{align}
    \max \; \sum_{v\in\V} \min_{\l\in\L_v} \theta_v'^\phi(\l)\span \\
    \forall uv\in\E,\,(\k,\l)\in\L_u\times\L_v:\;& \theta^\phi_{uv}(\k,\l)\geq0\\
    \forall v\in\V,\, u\in \N_v,\,\l\in\L_v:\;&\phi_{v\rightarrow u}(\l)\in\mathbb{R}.
\end{align}
\end{subequations}
Even though LP relaxation in this form was considered, e.g., in~\cite{cooper2007optimal}, computer-vision literature~\cite{hutschenreiter2021fusion,zhang2016pairwise,werner2007linear,savchynskyy2019discrete,tourani2018mplp++} typically considers the unconstrained form
\begin{subequations}\label{eq:dual_unconstrained_MRF_LP_relax}
\begin{align}
    \max \; \sum_{v\in\V} \min_{\l\in\L_v} \theta_v'^\phi(\l) + \sum_{uv\in\E} \min_{(\k,\l)\in\L_u\times\L_v}\theta^\phi_{uv}(\k,\l)\span \\
    \forall v\in\V,\, u\in \N_v,\,\l\in\L_v:\;&\phi_{v\rightarrow u}(\l)\in\mathbb{R}.
\end{align}
\end{subequations}
It is known~\cite[Remark~3]{werner2007linear} that the optimal values of~\eqref{eq:dual_MRF_LP_relax} and~\eqref{eq:dual_unconstrained_MRF_LP_relax} coincide. In detail, one can introduce the constraint $\sum_{(\k,\l)\in\L_u\times\L_v}\mu_{uv}(\k,\l)=1$ for each~$uv\in\E$ into the primal~\eqref{eq:MRF_LP_relax} and eliminate the corresponding dual variables analogously to~\eqref{eq:optimal_alpha}. Note, such a change does not influence the optimal value because the additional constraints are already implied by primal constraints~\eqref{eq:MRF_LP_relax:d} and~\eqref{eq:MRF_LP_relax:e}.

\end{remark}

\section{Compared Methods and Experimental Evaluation}\label{se:compared_methods}

Given the background from the previous sections, let us now summarize the algorithms that will be compared. A general algorithmic scheme is shown in Algorithm~\ref{al:general_scheme}. There, the dual variables~$\beta$ and~$\phi$ are first initialized to be feasible (by assigning them zero values). Then, we iteratively improve the current dual solution by separately updating the~$\phi$ and~$\beta$ variables. We do not consider the~$\alpha$ variables since their value is assumed to be implicitly determined by~\eqref{eq:optimal_alpha} based on the current values of~$\beta$ and~$\phi$.

Since the common optimal value of the primal-dual pair~\eqref{eq:iQAP_LP_relax} constitutes a lower bound on the optimal value of the original IQAP problem~\eqref{eq:iQAP}, any dual-feasible solution also provides a lower bound, which is determined by its objective. Algorithm~\ref{al:general_scheme} can be thus seen as gradually improving a lower bound on the optimal value of the IQAP problem.

\subsection{Compared Algorithms}\label{se:compared_algs}
To obtain a concrete algorithm, we need to define how the individual steps of Algorithm~\ref{al:general_scheme} are performed precisely:
\begin{itemize}
    \item Concerning line~3, the $\phi$ variables are improved by one loop of MPLP algorithm in~\cite{zhang2016pairwise} and by one loop of MPLP++ in~\cite{hutschenreiter2021fusion}. Since MPLP++~\cite{tourani2018mplp++} is an improved\footnote{ Although the fixed points of MPLP and MPLP++ coincide, it was shown that, if initialized at the same point, the dual objective after a single MPLP++ iteration is not worse than the dual objective after a single MPLP iteration~\cite[Section~5.1]{tourani2018mplp++}.} version of MPLP~\cite{globerson2007fixing}, we consider only MPLP++ in our experiments, i.e., line~3 is performed by running a single loop of MPLP++ algorithm.

    \item The update of the other dual variables on line~4 can be again done approximately by a single loop of BCA updates (as in~\cite{hutschenreiter2021fusion}) or exactly, e.g., by the Hungarian method (as in~\cite{zhang2016pairwise}). In addition to these options, we consider relative-interior updates, which can be obtained by combining any exact method with our Algorithm~\ref{al:ri-solution}.
\end{itemize}

\begin{algorithm}[t]\SetAlgoVlined\SetKwInOut{Output}{output} 
\Input{instance of the IQAP}
\Output{lower bound on the optimal value}
Initialize~$\beta:=0$ and~$\phi:=0$.\\
\Repeat{\normalfont{termination condition is met}}{
Improve~$\phi$ variables for fixed~$\beta$ variables.\\
Improve~$\beta$ variables for fixed~$\phi$ variables.
}
\Return{\normalfont{current dual objective~\eqref{eq:iQAP_LP_relax:a} where~$\alpha$ is~\eqref{eq:optimal_alpha}}}
\caption{General algorithmic scheme for obtaining a lower bound on the IQAP.}
\label{al:general_scheme}
\end{algorithm}

As listed above, the only difference in the compared methods lies in how the update of the dual variables~$\beta$ is performed. To be precise, the three options for updating~$\beta$ variables yield the following algorithms:
\begin{itemize}
    \item \bcaalg/~performs a single loop of coordinate-wise updates of~$\beta_\l$ while adhering to the relative-interior rule~\cite{werner2020relative}. In detail, we sequentially perform the updates~\eqref{eq:coordinate_ascent} for all~$\l\in\L\setminus\{\#\}$. Note that, after a loop of updates\footnote{It is likely that even if the updates were performed repeatedly, the objective need not converge to the optimum. This follows from the fact that the dual of LP formulation of LAP is not solvable by BCA~\cite[Section~5.1]{dlask2022classes} (although the optimum was frequently reached in experiments), so the dual of LP formulation of ILAP is also likely not solvable by BCA.}, the resulting~$\beta$ need not be optimal for the ILAP subproblem~\eqref{eq:iLAP_cpaf}.
    \item \hungalg/~updates the dual variables~$\beta$ by exactly solving the ILAP subproblem. For this, we use the reduction to the LAP from Section~\ref{se:reduction}, solve the resulting LAP by the Hungarian method~\cite{kuhn1955hungarian}, and obtain an optimal solution of the dual LP formulation of the ILAP via~\eqref{eq:alpha_beta}. To make the iterations faster, we store the optimal values of $\alpha$ and $\beta$ variables of the LAP subproblem and, in the next iteration, we initialize these variables in the Hungarian method to the values from the previous iteration and, if necessary, decrease their values to make them feasible.
    \item \hungrialg/~is the same as \hungalg/~but when the dual solution of the LAP is obtained, we shift it to the relative interior of optimizers using Algorithm~\ref{al:ri-solution}. By Theorem~\ref{th:maps_ri_to_ri}, this choice of optimizer adheres to the relative-interior rule.
\end{itemize}

For completeness, we note that both MPLP++~\cite{tourani2018mplp++} and MPLP~\cite{globerson2007fixing} do not adhere to the relative-interior rule (see \cite{werner2020relative}). There exist other methods that perform updates of the individual $\phi$ variables to the relative interior of optimizers, e.g. max-sum diffusion, but these options require longer time to converge in practice \cite{werner2020relative,tourani2018mplp++}. All of these methods attain the same fixed points when applied to the dual~\eqref{eq:MRF_LP_relax}, however, these fixed points are not guaranteed to be optimal for~\eqref{eq:MRF_LP_relax}, so the values of the $\phi$ variables are also not guaranteed to be optimal already after a single update. Consequently, the update on line~3 of Algorithm~\ref{al:general_scheme} performs only several individual BCA updates instead of finding the optimal values of $\phi$ for the fixed values of~$\beta$.

Still, the bounds computed by Algorithm~\ref{al:general_scheme} may in general depend on how the values on line~3 are chosen as the dual~\eqref{eq:MRF_LP_relax} is only a subproblem of~\eqref{eq:iQAP_LP_relax}. However, we restrict our focus only on the differences caused by choosing different methods for solving the ILAP subproblem~\eqref{eq:iLAP_LP_relax_subproblem}.

\begin{remark}
{\normalfont \bcaalg/} is analogous to the BCA algorithm considered in~\cite{hutschenreiter2021fusion} and {\normalfont \hungalg/} is similar to Hungarian-BP~\cite{zhang2016pairwise} except that we use MPLP++ instead of MPLP. Another difference to the particular implementations of these methods lies in the fact that the precise form of the duals in~\cite{hutschenreiter2021fusion} and~\cite{zhang2016pairwise} is different (and also different from the dual that we consider here, which is more compact). In more detail, our dual can be obtained from the ones in~\cite{hutschenreiter2021fusion} or~\cite{zhang2016pairwise} by adding more constraints and eliminating variables. It is known~\cite{dlask2022classes} that the quality of fixed points of BCA highly depends on the precise problem formulation. In the aforementioned case, this change of formulation seems to have a negative impact on the quality of fixed points.
\end{remark}

\subsection{Problem Instances Used for Evaluation}
For evaluation, we used the recent computer-vision benchmark~\cite{haller2022comparative} (451~instances), as well as the operations-research benchmark QAPLIB~\cite{burkard1997qaplib} (132~instances). This resulted in 583~instances in total.\footnote{We downloaded the computer-vision instances in \texttt{.dd} format from \url{https://vislearn.github.io/gmbench/datasets/}.
The datasets considered in~\cite{hutschenreiter2021fusion} constitute a strict subset of those in~\cite{haller2022comparative}. Namely, the groups \instancegroup{caltech-small}, \instancegroup{caltech-large}, and \instancegroup{house-sparse} are missing in~\cite{hutschenreiter2021fusion}. We converted the QAPLIB benchmark to \texttt{.dd} format and subtracted large constants from the unary costs to guarantee a complete assignment when solving the IQAP. In our experimental evaluation, we add the constants back to the bound values so that they are comparable with other sources in the literature.} The sizes and densities of the instances from different groups are shown in Table~\ref{ta:instances} where the first 11 groups contain computer-vision instances and the remaining 15 groups are instances from QAPLIB. Due to the different nature of the instances, we separate them visually in our tables.

\begin{table}[t]
    \centering
    \begin{tabular}{ccccccc}
    \multirow{2}{*}{group}& \multirow{2}{*}{instances} &\multirow{2}{*}{$|\V|$} &\multirow{2}{*}{$|\L\setminus\{\#\}|$} &\multirow{2}{*}{$|\L_v\setminus\{\#\}|$}&\multirow{2}{*}{$|\E|$}&density~$\tfrac{2|\E|}{|\V|(|\V|-1)}$ \\ 
    &  &  &   & && (in percent) \\ \hline  \hline
caltech-large & 9       &  36-219 & 51-341 & 1-62 & 612-7448 & 31.2-97.1 \\ \hline  
caltech-small & 21      &  9-117 & 14-201 & 1-67 & 36-2723 & 30.6-100  \\ \hline
car & 30      &  19-49 & 19-49 & 19-49 & 46-131 & 11.1-26.9 \\ \hline  
flow & 6       &  48-126 & 51-130 & 1-19 & 1100-5227 & 44.6-97.5 \\ \hline  
hotel & 105      &  30 & 30 & 30 & 435 & 100 \\ \hline  
house-dense & 105      &  30 & 30 & 30 & 435 & 100 \\ \hline  
house-sparse & 105       &  30 & 30 & 30 & 79 & 18.2  \\ \hline
motor & 20       &  15-52 & 15-52 & 15-52 & 33-139 & 10.5-32.4 \\ \hline
opengm & 4       &  19-20 & 19-20 & 19-20 & 112-190 & 65.5-100 \\ \hline  
pairs & 16       &  511-565 & 523-565 & 20-24 & 25173-34334 & 18.2-22.8 \\ \hline  
worms & 30       &  558 & 1202-1427 & 20-127 & 2343-2363 & 1.51-1.52 \\ \hline\hline  
bur & 8       &  26 & 26 & 26 & 325 & 100 \\ \hline  
chr & 14      &  12-25 & 12-25 & 12-25 & 11-24 & 8-16.7 \\ \hline  
els & 1       &  19 & 19 & 19 & 171 & 100 \\ \hline  
esc & 18       &  16-64 & 16-64 & 16-64 & 0-141 & 0-95.8 \\ \hline  
had & 5       &  12-20 & 12-20 & 12-20 & 66-190 & 100 \\ \hline  
kra & 3       &  30-32 & 30-32 & 30-32 & 165-435 & 33.3-100 \\ \hline  
lipa & 16       &  20-90 & 20-90 & 20-90 & 179-4005 & 94.2-100 \\ \hline  
nug & 15       &  12-30 & 12-30 & 12-30 & 66-435 & 66.4-100 \\ \hline  
rou & 3       &  12-20 & 12-20 & 12-20 & 66-189 & 99.5-100 \\ \hline  
scr & 3       &  12-20 & 12-20 & 12-20 & 28-62 & 32.6-42.4 \\ \hline  
sko & 13       &  42-100 & 42-100 & 42-100 & 861-4950 & 100 \\ \hline  
ste & 3       &  36 & 36 & 36 & 630 & 100 \\ \hline  
tai & 26       &  10-100 & 10-100 & 10-100 & 45-4950 & 3.87-100 \\ \hline  
tho & 2       &  30-40 & 30-40 & 30-40 & 435-780 & 100 \\ \hline  
wil & 2      &  50-100 & 50-100 & 50-100 & 1225-4950 & 100  \\ \hline
    \end{tabular}
    \caption{Properties of instances.}
    \label{ta:instances}
\end{table}

To possibly avoid trivial fixed points, we adjust each instance in the following way: the value~$\theta_{uv}(\l,\l)=0$ is replaced by~$10^7$ for all~$uv\in\E$ and~$\l\in(\L_v\cap \L_u) \setminus \{\#\}$, which was also done in~\cite{hutschenreiter2021fusion}. This change does not influence the optimal value of the IQAP instance because, for any~$x$ feasible to~\eqref{eq:iQAP} and any~$uv\in \E$, $x_u \in\L_u\setminus\{\#\}$ implies~$x_u\neq x_v$. On the other hand, this may increase the optimal value of the LP relaxation~\eqref{eq:iQAP_LP_relax}, i.e., improve the bound that it provides.

\begin{figure}
 \centering
  \includegraphics[width=0.7\textwidth]{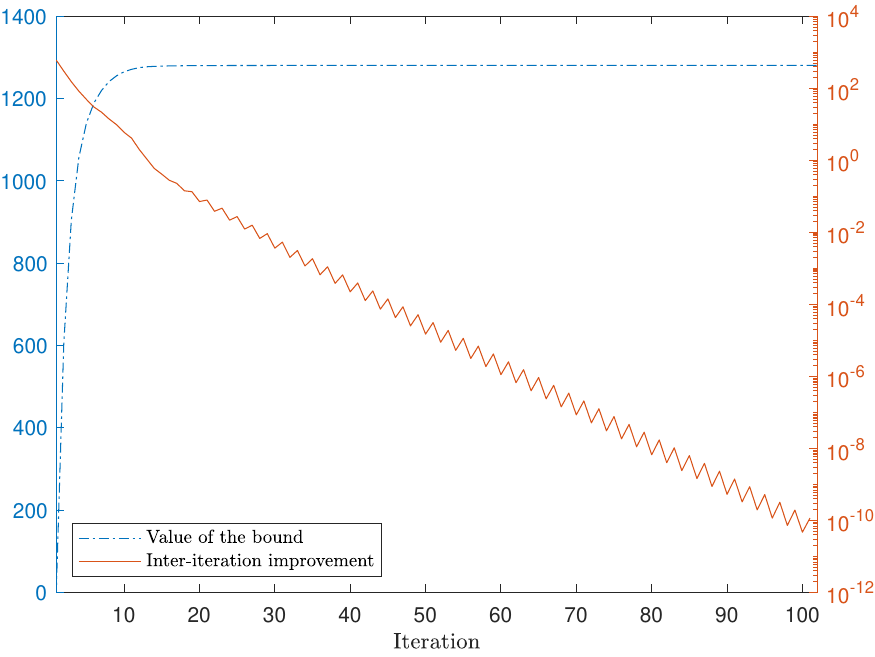}
 \caption{Current bound and relative inter-iteration improvement of the bound (the latter in logarithmic scale) for the instance \textit{lipa20a}. Only the first 100 iterations of $\hungalg/$ are shown.
 }\label{fig:convergence}
 \end{figure}
 
\subsection{Stopping Criteria}\label{sec:stopping-criteria}
Let us emphasise that the iterative scheme in Algorithm~\ref{al:general_scheme} need not attain a fixed point after a finite number of iterations. The obtained results therefore depend on the chosen stopping condition (see line~5 of the algorithm), which allows to balance between different trade-offs between runtime and quality of the bound. We exemplify such progress on the instance \textit{lipa20a} in Figure~\ref{fig:convergence} which shows the current bound attained by $\hungalg/$ for each of the first 100 iterations and also the inter-iteration improvements of the bound.

For the purpose of comparison of the aforementioned methods with an LP solver in Section~\ref{se:comparison_with_LP_solver}, we use the following stopping condition. With $b_i$ denoting the bound attained by the algorithm in iteration~$i$, we terminate the algorithm if, for 10 consecutive iterations, $b_i-b_{\lceil i/2 \rceil} \leq p (b_{\lceil i/2 \rceil}-b_1)$ where $p>0$ is a small fixed constant. In words, we should terminate if the improvement in the last half of the iterations is less than $p$ times the overall improvement in the first half of the iterations. In our case, we use $p\in\{0.01,0.1\}$. For example, in Figure~\ref{fig:convergence}, this results in 16 and 28 iterations, respectively. Finally, note that if the number of iterations (and thus also the runtime) is halved, then the obtained bound gets worse only by a factor of $p$, i.e., by 1\% and 10\% in our case.

To compare the quality of the bounds attainable by the methods in Section~\ref{se:compare_attainable_bounds}, we consider a criterion based on absolute improvements, i.e., terminate if $b_{i+1}-b_{i}\leq q$ for some small value of~$q>0$ for 10 consecutive iterations. We do not base our criteria on relative improvements (e.g., $\frac{b_{i+1}-b_{i}}{b_i}\leq q'$) because the bounds can be positive, negative, or zero.

Finally, we compare the performance of \hungrialg/ and \hungalg/ in Section~\ref{se:detailed_comparison} by comparing the bound computed after a common time limit or comparing the runtimes necessary to achieve a target bound.

\subsection{Comparison with LP Solver}\label{se:comparison_with_LP_solver}
We implemented the corresponding versions of Algorithm~\ref{al:general_scheme} in Matlab and compared them with Gurobi 11.0.1~\cite{gurobi}, which is a general-purpose LP solver. The evaluation was performed on a laptop with i7-7500U processor at 2.7~GHz and 16~GB RAM. In 47 out of the 583 instances, Gurobi ran out of memory or did not manage to compute the optimal value within the pre-specified time limit of 90 minutes per instance.\footnote{In some cases, the out-of-memory error happened only after the presolve phase of the solver finished. As 90~minutes did not suffice even for the presolve phase in some instances, we do not know in how many instances the memory would not be sufficient if longer runtime was allowed. For this reason, we do not distinguish the reason for not being solved by the LP solver.}

To compare the bounds computed by the considered methods with the optimal value of the linear program~\eqref{eq:iQAP_LP_relax}, we denote by~$B_\bcaalg/$, $B_\hungalg/$, $B_\hungrialg/$, and $B_\gurobialg/$ the bound attained by each method on a particular instance and the bound equal to the optimal value of the linear program (computed by \gurobialg/), respectively. The normalized bound of a method $m\in\{\bcaalg/, \hungalg/,\hungrialg/\}$ on the instance is then defined as
\begin{equation}\label{eq:normalized_bound}
    \frac{B_m-B_{\texttt{init}}}{B_\gurobialg/-B_{\texttt{init}}}
\end{equation}
where
\begin{equation}\label{eq:initial_bound}
    B_{\texttt{init}} = \sum_{v\in\V} \min_{\l\in\L_v} \theta_v(\l) + \sum_{uv\in\E} \min_{(\k,\l)\in\L_u\times\L_v}\theta_{uv}(\k,\l)
\end{equation}
is the initial bound from which the methods $\bcaalg/$, $\hungalg/$, and $\hungrialg/$ start. In other words, \eqref{eq:normalized_bound} scales the bounds so that they become comparable across instances. For example, value 1 indicates that the method attained optimum\footnote{For 8 instances, it happened that $B_\gurobialg/<0$ despite $B_{\texttt{init}}=0$, which can be expected as \gurobialg/ computes the result only up to a certain precision. The lowest of these values is $-8.5\cdot 10^{-4}$, so we compare the computed bounds up to a threshold of $10^{-3}$. I.e., we set the normalized bound to 1 if $B_\gurobialg/=B_{\texttt{init}}$ (up to the threshold $10^{-3}$), which happened for 71 instances. Also, the normalized bound is equal to 1 if the difference $B_\gurobialg/-B_m$ is less than $10^{-3}$.} of the linear program whereas value 0 indicates that the method was not able to improve upon the trivial bound~\eqref{eq:initial_bound}. In case that $\gurobialg/$ could not solve the particular instance, we replace $B_\gurobialg/$ by $\max\{B_\bcaalg/, B_\hungalg/,B_\hungrialg/\}$ in \eqref{eq:normalized_bound}.

We report the group-wise average normalized bounds in Tables~\ref{ta:stopping_1p} and \ref{ta:stopping_10p} for the stopping criterion with $p=0.01$ and $p=0.1$, respectively. The groups \instancegroup{flow}, \instancegroup{lipa}, \instancegroup{sko}, \instancegroup{tai}, and \instancegroup{wil} are further divided in each of the tables -- here, the first line corresponds to the instances that were solved by $\gurobialg/$ and the second line to the remaining ones. In addition, \gurobialg/ could not solve any instance from the group \instancegroup{pairs}. We do not state the normalized bound for \gurobialg/ as it would be either equal to 1 or undefined, depending on whether \gurobialg/ managed to solve the instance.

Regarding the attained normalized bounds, $\hungrialg/$ and $\hungalg/$ typically attain similar bounds which are often superior to (or at least comparable to) $\bcaalg/$. Despite being similar, $\hungrialg/$ attains better average normalized bounds than $\hungalg/$  when $p=0.01$, with the only exception being groups \instancegroup{house-dense} and \instancegroup{chr} and 12 cases where the averages are identical. For the more strict stopping criterion with $p=0.1$, $\hungalg/$ improves upon $\hungrialg/$ also in \instancegroup{hotel}, \instancegroup{motor}, \instancegroup{pairs}, and \instancegroup{tai}. Independently of the stopping criterion, \bcaalg/ is worse than the other methods and, e.g., cannot improve the initial bound in any instance within the groups \instancegroup{bur} and \instancegroup{rou}.

\newcommand\ExtraSep
{\dimexpr\cmidrulewidth+\aboverulesep+\belowrulesep\relax}

\afterpage{
\begin{landscape}

\begin{table}[t]
    \centering
    \resizebox{1.4\textwidth}{!}{\begin{tabular}{cc|PPP|cccc|cccc}
    \multirow{2}{*}{group}& \multirow{2}{*}{\begin{tabular}{c}instances solved \\ by \gurobialg/ \end{tabular}} & \multicolumn{3}{c|}{average normalized bound} & \multicolumn{4}{c|}{arithmetic mean of runtimes (seconds)} & \multicolumn{4}{c}{\begin{tabular}{c}shifted geometric mean \\
                of runtimes (seconds)
                \end{tabular}}\\ \cline{3-13}
    &  & \bcaalg/ & \hungalg/ & \hungrialg/  & \bcaalg/ & \hungalg/ & \hungrialg/ & \gurobialg/ & \bcaalg/ & \hungalg/ & \hungrialg/ & \gurobialg/ \\ \hline \hline
caltech-large & 9/9& 0.8704 & 0.9461 & \textbf{0.9476} & 14.79 & \textbf{12.43} & 12.68 & 1873.91 & 13.05 & \textbf{10.53} & 10.94 & 1486.74 \\ \hline 
caltech-small & 21/21& 0.8898 & 0.9541 & \textbf{0.9555} & 3.85 & \textbf{2.78} & 2.97 & 294.27 & 3.51 & \textbf{2.56} & 2.73 & 64.60 \\ \hline 
car & 30/30& 0.9825 & 0.9864 & \textbf{0.9882} & 4.85 & \textbf{0.61} & 0.78 & 4.78 & 4.53 & \textbf{0.60} & 0.77 & 4.38 \\ \hline 
\multirow{2}{*}{flow} & 3/6& \textbf{0.9833} & 0.9817 & 0.9820 & 17.06 & \textbf{13.54} & 14.43 & 983.67 & 15.96 & \textbf{11.55} & 12.01 & 462.13 \\ 
& 0/6& 0.9419 & 0.9881 & \textbf{0.9985} & 68.26 & \textbf{56.24} & 97.20 & -- & 63.99 & \textbf{55.73} & 89.32 & -- \\ \hline 
hotel & 105/105& 0.9869 & 0.9940 & \textbf{0.9945} & 4.14 & \textbf{1.63} & 1.99 & 145.39 & 3.40 & \textbf{1.46} & 1.69 & 140.50 \\ \hline 
house-dense & 105/105& 0.9996 & \textbf{0.9998} & 0.9993 & 7.21 & \textbf{3.14} & 3.54 & 173.00 & 7.01 & \textbf{3.08} & 3.47 & 172.83 \\ \hline 
house-sparse & 105/105& \textbf{1.0000} & \textbf{1.0000} & \textbf{1.0000} & 1.56 & \textbf{0.25} & 0.32 & 2.86 & 1.56 & \textbf{0.25} & 0.32 & 2.86 \\ \hline 
motor & 20/20& 0.9882 & 0.9919 & \textbf{0.9922} & 2.73 & \textbf{0.39} & 0.49 & 3.75 & 2.50 & \textbf{0.39} & 0.48 & 3.37 \\ \hline 
opengm & 4/4& 0.8590 & 0.8743 & \textbf{0.8746} & 1.05 & \textbf{0.56} & 0.63 & 5.74 & 1.04 & \textbf{0.56} & 0.63 & 5.52 \\ \hline 
pairs & 0/16& 0.7851 & 0.9990 & \textbf{0.9993} & \textbf{523.10} & 1349.60 & 1271.49 & -- & \textbf{508.14} & 1343.42 & 1262.57 & -- \\ \hline 
worms & 30/30& 0.9962 & 0.9957 & \textbf{0.9966} & \textbf{26.82} & 519.44 & 424.81 & 451.18 & \textbf{26.49} & 492.07 & 412.47 & 450.52 \\ \hline \hline 
bur & 8/8& 0.0000 & 0.7109 & \textbf{0.7115} & \textbf{0.90} & 2.72 & 3.02 & 31.04 & \textbf{0.90} & 2.71 & 3.01 & 28.28 \\ \hline 
chr & 14/14& 0.4196 & \textbf{0.9665} & 0.9661 & \textbf{0.22} & \textbf{0.22} & 0.26 & 0.77 & \textbf{0.22} & \textbf{0.22} & 0.26 & 0.77 \\ \hline 
els & 1/1& \textbf{0.0000} & \textbf{0.0000} & \textbf{0.0000} & 0.37 & \textbf{0.17} & \textbf{0.17} & 4.13 & 0.37 & \textbf{0.17} & \textbf{0.17} & 4.13 \\ \hline 
esc & 18/18& \textbf{1.0000} & \textbf{1.0000} & \textbf{1.0000} & 1.03 & \textbf{0.10} & 0.12 & 1.96 & 0.91 & \textbf{0.10} & 0.12 & 1.80 \\ \hline 
had & 5/5& 0.5098 & 0.8962 & \textbf{0.8965} & \textbf{0.26} & 0.48 & 0.52 & 3.45 & \textbf{0.26} & 0.48 & 0.51 & 3.24 \\ \hline 
kra & 3/3& \textbf{1.0000} & \textbf{1.0000} & \textbf{1.0000} & 1.20 & \textbf{0.41} & 0.45 & 20.03 & 1.20 & \textbf{0.41} & 0.45 & 18.00 \\ \hline 
\multirow{2}{*}{lipa} & 6/16& 0.4269 & 0.9088 & \textbf{0.9102} & \textbf{1.90} & 2.71 & 2.73 & 630.65 & \textbf{1.80} & 2.58 & 2.63 & 153.93 \\  
 & 0/16& 0.2565 & 0.9983 & \textbf{0.9992} & \textbf{39.89} & 296.90 & 293.90 & -- & \textbf{30.73} & 151.53 & 150.16 & -- \\ \hline 
nug & 15/15& \textbf{0.9986} & \textbf{0.9986} & \textbf{0.9986} & 0.52 & \textbf{0.20} & 0.23 & 6.91 & 0.51 & \textbf{0.20} & 0.23 & 5.92 \\ \hline 
rou & 3/3& 0.0000 & 0.8186 & \textbf{0.8199} & \textbf{0.25} & 0.64 & 0.62 & 3.32 & \textbf{0.25} & 0.64 & 0.62 & 3.08 \\ \hline 
scr & 3/3& \textbf{0.9947} & \textbf{0.9947} & \textbf{0.9947} & 0.19 & \textbf{0.05} & 0.10 & 0.98 & 0.19 & \textbf{0.05} & 0.10 & 0.98 \\ \hline 
\multirow{2}{*}{sko} & 4/13& \textbf{1.0000} & \textbf{1.0000} & \textbf{1.0000} & 7.80 & \textbf{3.44} & 3.64 & 1114.86 & 7.45 & \textbf{3.32} & 3.48 & 749.46 \\  
 & 0/13& \textbf{1.0000} & \textbf{1.0000} & \textbf{1.0000} & 70.35 & 65.77 & \textbf{53.45} & -- & 64.91 & 55.56 & \textbf{47.51} & -- \\ \hline 
ste & 3/3& \textbf{1.0000} & \textbf{1.0000} & \textbf{1.0000} & 2.20 & \textbf{0.83} & 0.85 & 86.66 & 2.20 & \textbf{0.83} & 0.85 & 86.59 \\ \hline 
\multirow{2}{*}{tai} & 18/26& 0.0556 & 0.5519 & \textbf{0.5543} & \textbf{1.33} & 1.56 & 2.07 & 194.73 & \textbf{1.21} & 1.41 & 1.86 & 36.70 \\  
 & 0/26& 0.5000 & 0.9997 & \textbf{0.9998} & \textbf{49.03} & 560.05 & 569.40 & -- & \textbf{30.71} & 80.18 & 90.99 & -- \\ \hline 
tho & 2/2& \textbf{1.0000} & \textbf{1.0000} & \textbf{1.0000} & 2.23 & 0.86 & \textbf{0.84} & 102.01 & 2.21 & 0.85 & \textbf{0.84} & 90.78 \\ \hline 
\multirow{2}{*}{wil} & 1/2& \textbf{1.0000} & \textbf{1.0000} & \textbf{1.0000} & 6.93 & 9.60 & \textbf{6.69} & 611.07 & 6.93 & 9.60 & \textbf{6.69} & 611.07 \\  
 & 0/2& \textbf{1.0000} & \textbf{1.0000} & \textbf{1.0000} & 68.90 & 67.26 & \textbf{60.61} & -- & 68.90 & 67.26 & \textbf{60.61} & -- \\ \hline 
       \end{tabular}}
    \caption{Group-wise aggregated comparison of the methods in terms of the average normalized bound and shifted geometric means of runtime, case with $p=0.01$. See Section~\ref{se:comparison_with_LP_solver} for the details on the bound value representation and Section~\ref{sec:stopping-criteria} on the runtime. Groups not fully solved by Gurobi are split into two subgroups: solved and unsolved instances, as requiring different baselines for the bound normalization.}
    \label{ta:stopping_1p}
\end{table}

\begin{table}[t]
    \centering
    \resizebox{1.4\textwidth}{!}{\begin{tabular}{cc|PPP|cccc|cccc}
    \multirow{2}{*}{group}& \multirow{2}{*}{\begin{tabular}{c}instances solved \\ by \gurobialg/ \end{tabular}} & \multicolumn{3}{c|}{average normalized bound} & \multicolumn{4}{c|}{arithmetic mean of runtimes (seconds)} & \multicolumn{4}{c}{\begin{tabular}{c}shifted geometric mean \\
                of runtimes (seconds)
                \end{tabular}}\\ \cline{3-13}
    &  & \bcaalg/ & \hungalg/ & \hungrialg/  & \bcaalg/ & \hungalg/ & \hungrialg/ & \gurobialg/ & \bcaalg/ & \hungalg/ & \hungrialg/ & \gurobialg/ \\ \hline \hline
caltech-large & 9/9& 0.8474 & 0.9272 & \textbf{0.9280} & \textbf{3.64} & 4.16 & 3.90 & 1873.91 & \textbf{3.45} & 3.66 & 3.50 & 1486.74 \\ \hline 
caltech-small & 21/21& 0.8746 & 0.9420 & \textbf{0.9440} & 1.16 & \textbf{1.03} & 1.11 & 294.27 & 1.12 & \textbf{0.99} & 1.07 & 64.60 \\ \hline 
car & 30/30& 0.9658 & 0.9726 & \textbf{0.9731} & 1.54 & \textbf{0.20} & 0.25 & 4.78 & 1.50 & \textbf{0.20} & 0.25 & 4.38 \\ \hline 
\multirow{2}{*}{flow} & 3/6& 0.9337 & 0.9358 & \textbf{0.9369} & 2.32 & 2.11 & \textbf{2.07} & 983.67 & 2.30 & 2.07 & \textbf{2.04} & 462.13 \\  
 & 0/6& 0.9547 & 0.9970 & \textbf{1.0000} & \textbf{5.05} & 5.76 & 5.32 & -- & \textbf{4.99} & 5.53 & 5.19 & -- \\ \hline 
hotel & 105/105& 0.9706 & \textbf{0.9776} & 0.9775 & 1.39 & \textbf{0.62} & 0.65 & 145.39 & 1.39 & \textbf{0.62} & 0.65 & 140.50 \\ \hline 
house-dense & 105/105& 0.9170 & \textbf{0.9178} & 0.9171 & 1.36 & \textbf{0.59} & 0.62 & 173.00 & 1.36 & \textbf{0.59} & 0.62 & 172.83 \\ \hline 
house-sparse & 105/105& \textbf{1.0000} & \textbf{1.0000} & \textbf{1.0000} & 1.28 & \textbf{0.20} & 0.26 & 2.86 & 1.28 & \textbf{0.20} & 0.26 & 2.86 \\ \hline 
motor & 20/20& 0.9822 & \textbf{0.9856} & 0.9853 & 1.27 & \textbf{0.17} & 0.21 & 3.75 & 1.22 & \textbf{0.17} & 0.21 & 3.37 \\ \hline 
opengm & 4/4& 0.8505 & 0.8659 & \textbf{0.8666} & 0.38 & \textbf{0.22} & 0.24 & 5.74 & 0.37 & \textbf{0.22} & 0.24 & 5.52 \\ \hline 
pairs & 0/16& 0.7858 & \textbf{0.9996} & 0.9992 & \textbf{49.46} & 219.46 & 212.09 & -- & \textbf{49.25} & 217.63 & 210.30 & -- \\ \hline 
worms & 30/30& 0.9907 & 0.9905 & \textbf{0.9914} & \textbf{12.96} & 347.45 & 309.05 & 451.18 & \textbf{12.93} & 340.50 & 304.54 & 450.52 \\ \hline  \hline
bur & 8/8& 0.0000 & 0.6585 & \textbf{0.6590} & 0.90 & \textbf{0.60} & 0.65 & 31.04 & 0.90 & \textbf{0.60} & 0.65 & 28.28 \\ \hline 
chr & 14/14& 0.4196 & 0.9449 & \textbf{0.9457} & 0.22 & \textbf{0.08} & 0.10 & 0.77 & 0.22 & \textbf{0.08} & 0.10 & 0.77 \\ \hline 
els & 1/1& \textbf{0.0000} & \textbf{0.0000} & \textbf{0.0000} & 0.37 & \textbf{0.17} & \textbf{0.17} & 4.13 & 0.37 & \textbf{0.17} & \textbf{0.17} & 4.13 \\ \hline 
esc & 18/18& \textbf{1.0000} & \textbf{1.0000} & \textbf{1.0000} & 1.03 & \textbf{0.10} & 0.12 & 1.96 & 0.91 & \textbf{0.10} & 0.12 & 1.80 \\ \hline 
had & 5/5& 0.5098 & 0.8950 & \textbf{0.8956} & \textbf{0.26} & 0.28 & 0.31 & 3.45 & \textbf{0.26} & 0.28 & 0.31 & 3.24 \\ \hline 
kra & 3/3& \textbf{1.0000} & \textbf{1.0000} & \textbf{1.0000} & 1.20 & \textbf{0.41} & 0.45 & 20.03 & 1.20 & \textbf{0.41} & 0.45 & 18.00 \\ \hline 
\multirow{2}{*}{lipa} & 6/16& 0.4269 & 0.9063 & \textbf{0.9075} & 1.90 & 1.42 & \textbf{1.41} & 630.65 & 1.80 & 1.38 & \textbf{1.37} & 153.93 \\  
 & 0/16& 0.2589 & 0.9984 & \textbf{0.9993} & \textbf{39.89} & 105.79 & 105.02 & -- & \textbf{30.73} & 63.05 & 61.36 & -- \\ \hline 
nug & 15/15& \textbf{0.9986} & \textbf{0.9986} & \textbf{0.9986} & 0.52 & \textbf{0.20} & 0.23 & 6.91 & 0.51 & \textbf{0.20} & 0.23 & 5.92 \\ \hline 
rou & 3/3& 0.0000 & 0.8079 & \textbf{0.8131} & \textbf{0.25} & 0.29 & 0.33 & 3.32 & \textbf{0.25} & 0.29 & 0.33 & 3.08 \\ \hline 
scr & 3/3& \textbf{0.9947} & \textbf{0.9947} & \textbf{0.9947} & 0.19 & \textbf{0.05} & 0.10 & 0.98 & 0.19 & \textbf{0.05} & 0.10 & 0.98 \\ \hline 
\multirow{2}{*}{sko} & 4/13& \textbf{1.0000} & \textbf{1.0000} & \textbf{1.0000} & 7.80 & \textbf{3.44} & 3.64 & 1114.86 & 7.45 & \textbf{3.32} & 3.48 & 749.46 \\  
 & 0/13& \textbf{1.0000} & \textbf{1.0000} & \textbf{1.0000} & 70.35 & 65.77 & \textbf{53.45} & -- & 64.91 & 55.56 & \textbf{47.51} & -- \\ \hline 
ste & 3/3& \textbf{1.0000} & \textbf{1.0000} & \textbf{1.0000} & 2.20 & \textbf{0.83} & 0.85 & 86.66 & 2.20 & \textbf{0.83} & 0.85 & 86.59 \\ \hline 
\multirow{2}{*}{tai} & 18/26& 0.0556 & 0.5495 & \textbf{0.5514} & 1.33 & \textbf{0.90} & 1.33 & 194.73 & 1.21 & \textbf{0.86} & 1.23 & 36.70 \\ 
 & 0/26& 0.5000 & \textbf{0.9999} & 0.9997 & \textbf{49.03} & 229.39 & 245.51 & -- & \textbf{30.71} & 53.11 & 61.00 & -- \\ \hline 
tho & 2/2& \textbf{1.0000} & \textbf{1.0000} & \textbf{1.0000} & 2.23 & 0.86 & \textbf{0.84} & 102.01 & 2.21 & 0.85 & \textbf{0.84} & 90.78 \\ \hline 
\multirow{2}{*}{wil} & 1/2& \textbf{1.0000} & \textbf{1.0000} & \textbf{1.0000} & 6.93 & 9.60 & \textbf{6.69} & 611.07 & 6.93 & 9.60 & \textbf{6.69} & 611.07 \\  
 & 0/2& \textbf{1.0000} & \textbf{1.0000} & \textbf{1.0000} & 68.90 & 67.26 & \textbf{60.61} & -- & 68.90 & 67.26 & \textbf{60.61} & -- \\ \hline 
        \end{tabular}}
    \caption{Group-wise aggregated comparison of the methods in terms of the average normalized bound and shifted geometric means of runtime, case with $p=0.1$. See Section~\ref{se:comparison_with_LP_solver} for the details on the bound value representation and Section~\ref{sec:stopping-criteria} on the runtime. Groups not fully solved by Gurobi are split into two subgroups: solved and unsolved instances, as requiring different baselines for the bound normalization.}
    \label{ta:stopping_10p}
\end{table}

\end{landscape}}

Tables~\ref{ta:stopping_1p} and \ref{ta:stopping_10p} also report the arithmetic and shifted geometric means of runtimes within the groups. The latter is used as it is more robust against outliers \cite[Section~6.3.1]{beach2024enhancements}, we used the shift of 10~seconds, as \cite{beach2024enhancements}. The relative values of both aggregate measures of runtime are often similar within individual groups, with largest differences occurring in groups \instancegroup{caltech-small}, \instancegroup{flow}, \instancegroup{lipa}, \instancegroup{sko}, and \instancegroup{tai}.

The runtimes of \hungalg/ are generally comparable to those of \hungrialg/, with \hungrialg/ sometimes attaining even smaller mean runtime than \hungalg/ despite its computational overhead. Depending on the instance, this is caused by \hungrialg/ requiring less iterations or actually having lower runtime per iteration. We believe that the latter phenomenon occurs because, as described in Section~\ref{se:compared_algs}, we initialize the variables of the LAP subproblem to the values from the previous iteration (up to adjustments for feasibility) and the relative-interior solution may be more robust w.r.t.\ small changes of the problem and thus perhaps closer to the optimal values, making runtime of the Hungarian method lower. For example, in the group \instancegroup{worms} with $p=0.1$, \hungalg/ and \hungrialg/ take 20.3 seconds and 16.9 seconds per iteration on average, respectively.

The only case where \hungrialg/ is significantly slower than \hungalg/ are the 3 instances within the group \instancegroup{flow} with $p=0.01$ that were not solved by \gurobialg/. The reason is that \hungrialg/ attains a better bound and needs more iterations to reach it.

Generally, the number of iterations required by the individual methods varies. Among all the instances, \bcaalg/, \hungalg/, and \hungrialg/ need up to 16, 50, and 81 iterations with $p=0.1$ and up to 278, 222, and 496 iterations with $p=0.01$, respectively. The per-iteration runtime is within the intervals 4.0~ms -- 18~s, 1.6~ms -- 49~s, and 2.3~ms -- 39~s, respectively. Typically, \hungalg/ and \hungrialg/ need similar number of iterations, with \hungrialg/ having slightly more. \bcaalg/ requires less iterations than the other methods.

In groups \instancegroup{pairs}, \instancegroup{worms}, \instancegroup{bur}, \instancegroup{had}, \instancegroup{lipa}, \instancegroup{rou}, and \instancegroup{tai}, \bcaalg/ is significantly faster than \hungalg/ or \hungrialg/, but, with the only exception of \instancegroup{worms}, this is at the cost of much worse attained bounds. In all the other groups, the runtime of \bcaalg/ is higher than or similar to \hungalg/ and \hungrialg/.

Regarding the comparison of these methods to the LP solver, \hungalg/ and \hungrialg/ are typically orders of magnitude faster than \gurobialg/, except for the group \instancegroup{worms} where \bcaalg/ seems to be the best option. Note, all of these speedups are attained despite we used an unoptimized implementation in Matlab whereas \gurobialg/ is an efficient commercial solver implemented in C.

Beyond being faster, \hungalg/ and \hungrialg/ provide bounds comparable to those computed by \gurobialg/. E.g., with $p=0.01$, \bcaalg/, \hungalg/, and \hungrialg/ attain a normalized bound larger than 0.8 for 481, 515, and 516 instances, respectively, out of the 536 instances that \gurobialg/ solved. The normalized bound is larger than 0.99 for 410, 418, and 419 instances, respectively. For the case with $p=0.1$, these numbers decrease to 476, 512, 513 and 273, 274, 280, respectively.

\begin{figure}
\centering
\begin{subfigure}{0.7\textwidth}
\centering
    \includegraphics[width=\textwidth]{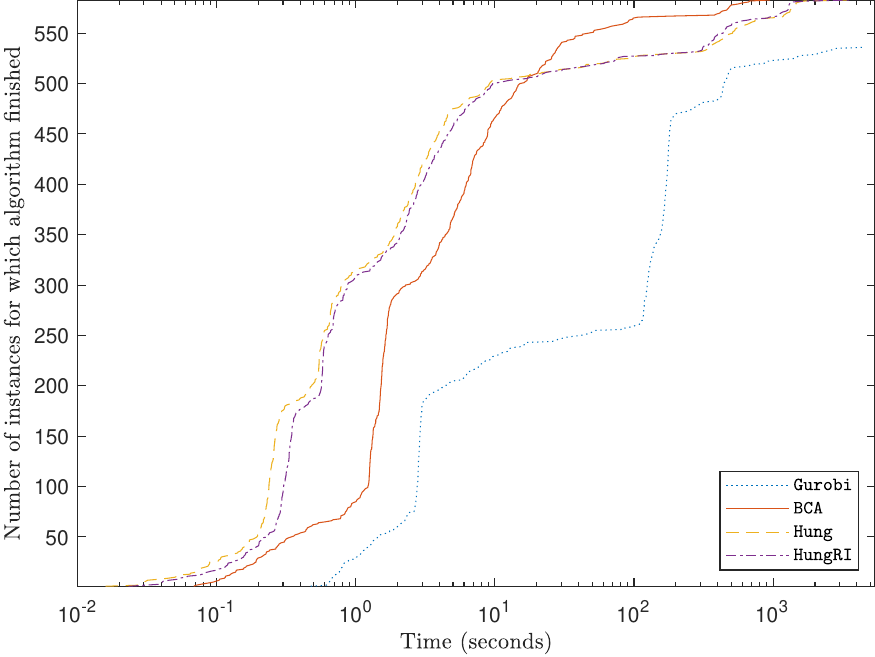}
    \caption{Case with $p=0.01$.}\label{fig:times_1percent}
\end{subfigure}

\vspace{1cm}

\begin{subfigure}{0.7\textwidth}
\centering
    \includegraphics[width=\textwidth]{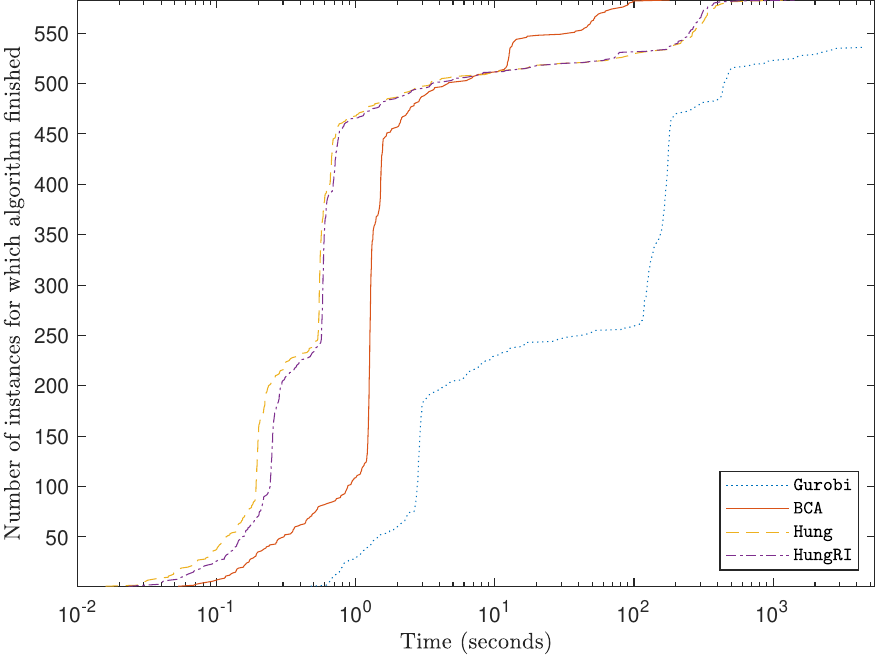}
     \caption{Case with $p=0.1$.}\label{fig:times_10percent}
\end{subfigure}
\caption{Cactus plots showing the distribution of the runtimes of the methods across all instances. See Section~\ref{sec:stopping-criteria} on the runtime.}
\label{fig:cactus}
\end{figure}

To provide additional comparison of the runtimes, we show in Figures~\ref{fig:times_1percent} and \ref{fig:times_10percent} the distribution of the runtimes of the individual methods for $p=0.01$ and $p=0.1$, respectively. Both figures confirm that \hungalg/ and \hungrialg/ have very similar runtimes. \bcaalg/ needs more time on the easier instances but is faster on the other end of the spectrum. \gurobialg/ is significantly slower than all the other methods.

\subsection{Comparison of the Quality of Fixed points}\label{se:compare_attainable_bounds}

In the previous section, we focused on stopping conditions that make the compared algorithms practical in terms of short overall runtime and high-quality bounds. In this section, we compare the quality of the attainable bounds of the methods to see how the choice of the relative interior point may improve the fixed points. For this, we let the methods run until the bound between the iterations improves by less than $10^{-7}$ in at least 10 consecutive iterations. For the computational time to be manageable, we increase the threshold to $1.5\cdot10^{-3}$ for the group \instancegroup{pairs}. Note, despite being theoretically incomparable, this condition typically turns out to be stronger than the one used in the previous section to allow for comparison of the (almost) converged bounds.

In Table~\ref{ta:bnd_stop_threshold}, we report, for each instance group and each method, how often it reached the best\footnote{As in Section~\ref{se:comparison_with_LP_solver}, a bound is considered best if it is not worse by more than $10^{-3}$ when compared to the actually best bound. In the group \instancegroup{pairs}, we use the same comparison threshold despite the weaker stopping condition. In this group, the bounds computed by \bcaalg/ are always worse than the best bounds by at least 375. The non-zero differences between the best bounds and the bounds computed by \hungrialg/ are at least 13.9 and the lowest non-zero differences between the best bounds and the bounds computed by \hungalg/ are 0.016, 1.42, and then larger than 4.20.} bound among the 3 methods and also the average bound for each group. The table shows that the average bound computed by \hungalg/ and \hungrialg/ is typically similar, but \hungrialg/ is most frequently capable of attaining the best bound.

Note that, for groups \instancegroup{els}, \instancegroup{esc}, \instancegroup{sko}, \instancegroup{ste}, \instancegroup{tho}, and \instancegroup{wil} within the QAPLIB benchmark, all methods return the zero bound for all instances. Except for the single instance in the group \instancegroup{els} (and perhaps the instances that could not be solved by \gurobialg/ and the optimal value is thus unknown to us), this value is optimal.

\afterpage{
\begin{landscape}
\begin{table}[t]
    \centering
    \begin{tabular}{cc|PPP|ccc}
    \multirow{2}{*}{group}& \multirow{2}{*}{\begin{tabular}{c}instances solved \\ by \gurobialg/ \end{tabular}} & \multicolumn{3}{c|}{\# best bound} & \multicolumn{3}{c}{average bound} \\ \cline{3-8}
    &  & \bcaalg/ & \hungalg/ & \hungrialg/  & \bcaalg/ & \hungalg/ & \hungrialg/ \\ \hline \hline
caltech-large & 9/9& 1 & 1 & \textbf{7} & -56348.01 & -47026.39 & \textbf{-46833.49} \\ \hline 
caltech-small & 21/21& 1 & 5 & \textbf{15} & -15770.24 & -12348.64 & \textbf{-12299.63} \\ \hline 
car & 30/30& 17 & 16 & \textbf{27} & -70.76 & -70.63 & \textbf{-70.56} \\ \hline 
\multirow{2}{*}{flow} & 3/6& \textbf{2} & \textbf{2} & \textbf{2} & -2420.64 & \textbf{-2418.64} & -2418.96 \\ 
 & 0/6& 0 & 1 & \textbf{2} & -3390.82 & -3360.56 & \textbf{-3357.64} \\ \hline 
hotel & 105/105& 91 & 98 & \textbf{102} & -4300.62 & -4298.20 & \textbf{-4298.08} \\ \hline 
house-dense & 105/105& 103 & \textbf{105} & 101 & -3778.24 & \textbf{-3778.19} & -3778.30 \\ \hline 
house-sparse & 105/105& \textbf{105} & \textbf{105} & \textbf{105} & \textbf{-66.78} & \textbf{-66.78} & \textbf{-66.78} \\ \hline 
motor & 20/20& 16 & \textbf{17} & \textbf{17} & -63.45 & \textbf{-63.34} & -63.35 \\ \hline 
opengm & 4/4& \textbf{2} & 1 & \textbf{2} & -177.35 & \textbf{-177.17} & \textbf{-177.17} \\ \hline 
pairs & 0/16& 0 & 3 & \textbf{13} & -70245.20 & -68368.80 & \textbf{-68366.13} \\ \hline 
worms & 30/30& \textbf{19} & 0 & 12 & -48500.45 & -48503.92 & \textbf{-48497.93} \\ \hline \hline
bur & 8/8& 0 & 1 & \textbf{7} & 0.00 & 253433.08 & \textbf{253447.92} \\ \hline 
chr & 14/14& 1 & 7 & \textbf{10} & 1726.43 & 6097.85 & \textbf{6112.43} \\ \hline 
els & 1/1& \textbf{1} & \textbf{1} & \textbf{1} & \textbf{0.00} & \textbf{0.00} & \textbf{0.00} \\ \hline 
esc & 18/18& \textbf{18} & \textbf{18} & \textbf{18} & \textbf{0.00} & \textbf{0.00} & \textbf{0.00} \\ \hline 
had & 5/5& 0 & 2 & \textbf{3} & 915.20 & 1556.29 & \textbf{1556.32} \\ \hline 
kra & 3/3& \textbf{3} & \textbf{3} & \textbf{3} & \textbf{12133.33} & \textbf{12133.33} & \textbf{12133.33} \\ \hline 
\multirow{2}{*}{lipa} & 6/16& 0 & 1 & \textbf{5} & 22820.83 & 43019.48 & \textbf{43044.85} \\ 
 & 0/16& 0 & 3 & \textbf{7} & 164006.60 & \textbf{618921.57} & 618792.19 \\ \hline 
nug & 15/15& \textbf{15} & \textbf{15} & \textbf{15} & \textbf{244.80} & \textbf{244.80} & \textbf{244.80} \\ \hline 
rou & 3/3& 0 & 0 & \textbf{3} & 0.00 & 54437.97 & \textbf{54523.46} \\ \hline 
scr & 3/3& \textbf{3} & \textbf{3} & \textbf{3} & \textbf{46761.33} & \textbf{46761.33} & \textbf{46761.33} \\ \hline 
\multirow{2}{*}{sko} & 4/13& \textbf{4} & \textbf{4} & \textbf{4} & \textbf{0.00} & \textbf{0.00} & \textbf{0.00} \\  
 & 0/13& \textbf{9} & \textbf{9} & \textbf{9} & \textbf{0.00} & \textbf{0.00} & \textbf{0.00} \\ \hline 
ste & 3/3& \textbf{3} & \textbf{3} & \textbf{3} & \textbf{0.00} & \textbf{0.00} & \textbf{0.00} \\ \hline 
\multirow{2}{*}{tai} & 18/26& 6 & 11 & \textbf{14} & 27083.33 & 71589.64 & \textbf{71654.02} \\  
 & 0/26& 4 & 5 & \textbf{7} & 0.00 & 143482.22 & \textbf{143538.70} \\ \hline 
tho & 2/2& \textbf{2} & \textbf{2} & \textbf{2} & \textbf{0.00} & \textbf{0.00} & \textbf{0.00} \\ \hline 
\multirow{2}{*}{wil} & 1/2& \textbf{1} & \textbf{1} & \textbf{1} & \textbf{0.00} & \textbf{0.00} & \textbf{0.00} \\  
 & 0/2& \textbf{1} & \textbf{1} & \textbf{1} & \textbf{0.00} & \textbf{0.00} & \textbf{0.00} \\ \hline 
    \end{tabular}
    \caption{Group-wise aggregated results for comparison of bounds attainable by our three methods with stronger stopping condition. See Section~\ref{se:compare_attainable_bounds} for details.}
    \label{ta:bnd_stop_threshold}
\end{table}
\end{landscape}}

\subsection{Detailed Comparison of \hungalg/ and \hungrialg/}\label{se:detailed_comparison}

We perform an additional detailed comparison of \hungalg/ and \hungrialg/ to better understand the potential trade-offs. Given the results from the previous sections, we compare the methods as follows:
\begin{itemize}
    \item For each instance, find out which method returned the worse (i.e., lower) bound and stop the other method when it reaches this bound. Then, for each group, compare which of the methods attained this bound faster and compute the group-wise shifted geometric mean of such runtimes. The results of this evaluation are in the left part of Tables~\ref{ta:additional_comparison_hungalgs2} and~\ref{ta:additional_comparison_hungalgs}.
    \item For each instance, find out which method terminated earlier and stop the other method at this time. Then, compare which of the methods reached better bound at this time and compute the group-wise averages of these bounds. The results of this evaluation are in the right part of Tables~\ref{ta:additional_comparison_hungalgs2} and~\ref{ta:additional_comparison_hungalgs}.
\end{itemize}

Tables~\ref{ta:additional_comparison_hungalgs2} and~\ref{ta:additional_comparison_hungalgs} are based on the stopping conditions from Section~\ref{se:comparison_with_LP_solver} with $p=0.1$ and Section~\ref{se:compare_attainable_bounds}, respectively. To better present the differences between the methods, let us note that we report the number of strictly faster terminations and strictly better bounds, which is different from the other tables in the paper.

Regarding Table~\ref{ta:additional_comparison_hungalgs2}, where the typically weaker stopping condition was used, the results of both methods are generally similar, sometimes with slight preferences towards one or the other method. E.g., \hungrialg/ seems preferable in the groups \instancegroup{pairs} and \instancegroup{worms} whereas \hungalg/ is preferable in \instancegroup{bur}.

In contrast, the differences between the methods get more significant in Table \ref{ta:additional_comparison_hungalgs}, where the stronger stopping condition was used. There, \hungrialg/ is more often capable of reaching the same bound \emph{significantly} faster than \hungalg/ -- this happened, e.g., in groups \instancegroup{caltech-large}, \instancegroup{caltech-small}, \instancegroup{worms}, and \instancegroup{rou}. In contrast, when \hungalg/ is faster, the gains in time are insignificant.

\afterpage{
\begin{landscape}

\begin{table}[t]
    \centering
    \begin{tabular}{cc|PP|Pc||PP|cc}
    \multirow{3}{*}{group}& \multirow{3}{*}{\# instances} & \multicolumn{4}{c||}{\begin{tabular}{c}target bound:  \\ stop the better method when it  \\ reaches the worse of the bounds \end{tabular}} & \multicolumn{4}{c}{\begin{tabular}{c}target time:\\ stop the slower method \\when the faster terminates\end{tabular}} \\ \cline{3-10}
    & & \multicolumn{2}{c|}{\# faster terminations} & \multicolumn{2}{c||}{\begin{tabular}{c}shifted geometric mean \\ of runtimes (seconds)
                \end{tabular}} & \multicolumn{2}{c|}{\# better bounds} & \multicolumn{2}{c}{average bound} \\ 
    &  & \hungalg/ & \hungrialg/ & \hungalg/ & \hungrialg/ & \hungalg/ & \hungrialg/ & \hungalg/ & \hungrialg/\\  \hline \hline
caltech-large & 9 & \textbf{5} & 4 & 3.66 & \textbf{3.44} & 3 & \textbf{6} & -49420.26 & \textbf{-49087.06} \\ \hline 
caltech-small & 21 & \textbf{16} & 5 & \textbf{0.99} & 1.02 & 9 & \textbf{12} & -12849.64 & \textbf{-12829.86} \\ \hline 
car & 30 & \textbf{29} & 1 & \textbf{0.20} & 0.24 & \textbf{23} & 4 & \textbf{-71.23} & -71.40 \\ \hline 
flow & 6 & \textbf{3} & \textbf{3} & 3.69 & \textbf{3.52} & 0 & \textbf{6} & -2942.15 & \textbf{-2939.35} \\ \hline 
hotel & 105 & \textbf{86} & 19 & \textbf{0.32} & 0.34 & \textbf{22} & 9 & \textbf{-4303.35} & -4303.44 \\ \hline 
house-dense & 105 & \textbf{97} & 8 & \textbf{0.59} & 0.62 & \textbf{66} & 39 & \textbf{-3803.16} & -3803.65 \\ \hline 
house-sparse & 105 & \textbf{102} & 3 & \textbf{0.12} & 0.15 & \textbf{2} & 0 & \textbf{-66.78} & \textbf{-66.78} \\ \hline 
motor & 20 & \textbf{16} & 4 & \textbf{0.15} & 0.17 & \textbf{11} & 1 & \textbf{-63.67} & -63.76 \\ \hline 
opengm & 4 & \textbf{4} & 0 & \textbf{0.22} & 0.24 & \textbf{3} & 1 & \textbf{-177.36} & -177.38 \\ \hline 
pairs & 16 & 1 & \textbf{15} & 217.09 & \textbf{210.30} & 4 & \textbf{12} & -69130.32 & \textbf{-69113.83} \\ \hline 
worms & 30 & 0 & \textbf{30} & 340.50 & \textbf{300.14} & 0 & \textbf{30} & -48585.75 & \textbf{-48555.80} \\ \hline \hline
bur & 8 & \textbf{6} & 2 & \textbf{0.60} & 0.65 & 2 & \textbf{6} & \textbf{237393.22} & 236908.74 \\ \hline 
chr & 14 & \textbf{13} & 1 & \textbf{0.07} & 0.09 & \textbf{11} & 1 & \textbf{5903.68} & 5861.14 \\ \hline 
els & 1 & 0 & \textbf{1} & \textbf{0.02} & \textbf{0.02} & \textbf{0} & \textbf{0} & \textbf{0.00} & \textbf{0.00} \\ \hline 
esc & 18 & \textbf{15} & 3 & \textbf{0.01} & \textbf{0.01} & \textbf{0} & \textbf{0} & \textbf{0.00} & \textbf{0.00} \\ \hline 
had & 5 & \textbf{3} & 2 & \textbf{0.26} & 0.29 & \textbf{3} & 2 & \textbf{1553.37} & 1552.84 \\ \hline 
kra & 3 & \textbf{3} & 0 & \textbf{0.04} & \textbf{0.04} & \textbf{0} & \textbf{0} & \textbf{12133.33} & \textbf{12133.33} \\ \hline 
lipa & 16 & 4 & \textbf{12} & 26.05 & \textbf{24.98} & 4 & \textbf{12} & 398904.14 & \textbf{399051.51} \\ \hline 
nug & 15 & \textbf{12} & 3 & \textbf{0.02} & \textbf{0.02} & \textbf{0} & \textbf{0} & \textbf{244.80} & \textbf{244.80} \\ \hline 
rou & 3 & 0 & \textbf{3} & 0.29 & \textbf{0.25} & 0 & \textbf{3} & 53820.10 & \textbf{53994.45} \\ \hline 
scr & 3 & \textbf{3} & 0 & \textbf{0.00} & 0.01 & \textbf{0} & \textbf{0} & \textbf{46761.33} & \textbf{46761.33} \\ \hline 
sko & 13 & 4 & \textbf{9} & 8.53 & \textbf{3.63} & \textbf{0} & \textbf{0} & \textbf{0.00} & \textbf{0.00} \\ \hline 
ste & 3 & \textbf{2} & 1 & 0.13 & \textbf{0.09} & \textbf{0} & \textbf{0} & \textbf{0.00} & \textbf{0.00} \\ \hline 
tai & 26 & \textbf{17} & 9 & 6.71 & \textbf{6.70} & 7 & \textbf{9} & \textbf{93055.79} & 93033.07 \\ \hline 
tho & 2 & \textbf{2} & 0 & \textbf{0.08} & 0.09 & \textbf{0} & \textbf{0} & \textbf{0.00} & \textbf{0.00} \\ \hline 
wil & 2 & 0 & \textbf{2} & 6.18 & \textbf{4.18} & \textbf{0} & \textbf{0} & \textbf{0.00} & \textbf{0.00} \\ \hline 
\end{tabular}
    \caption{Group-wise aggregated results of a detailed comparison of $\hungalg/$ and $\hungrialg/$ with different stopping conditions. The \emph{worse} bound and time has been computed based on the stopping condition described in Section~\ref{se:comparison_with_LP_solver} with $p=0.1$.
    See Section~\ref{se:detailed_comparison} for details. }
    \label{ta:additional_comparison_hungalgs2}
\end{table}

\begin{table}[t]
    \centering
    \begin{tabular}{cc|PP|Pc||PP|cc}
    \multirow{3}{*}{group}& \multirow{3}{*}{\# instances} & \multicolumn{4}{c||}{\begin{tabular}{c}target bound:  \\stop the better method when it  \\reaches the worse of the bounds \end{tabular}} & \multicolumn{4}{c}{\begin{tabular}{c}target time:\\stop the slower method \\when the faster terminates\end{tabular}} \\ \cline{3-10}
    & & \multicolumn{2}{c|}{\# faster terminations} & \multicolumn{2}{c||}{\begin{tabular}{c}shifted geometric mean \\ of runtimes (seconds)
                \end{tabular}} & \multicolumn{2}{c|}{\# better bounds} & \multicolumn{2}{c}{average bound} \\ 
    &  & \hungalg/ & \hungrialg/ & \hungalg/ & \hungrialg/ & \hungalg/ & \hungrialg/ & \hungalg/ & \hungrialg/ \\ \hline \hline 
caltech-large & 9 & 1 & \textbf{8} & 135.59 & \textbf{25.78} & 1 & \textbf{8} & -47026.40 & \textbf{-46833.49} \\ \hline 
caltech-small & 21 & 5 & \textbf{16} & 32.26 & \textbf{6.95} & 5 & \textbf{16} & -12349.31 & \textbf{-12299.63} \\ \hline 
car & 30 & 14 & \textbf{16} & 3.45 & \textbf{1.59} & 5 & \textbf{13} & -70.63 & \textbf{-70.57} \\ \hline 
flow & 6 & 2 & \textbf{4} & \textbf{65.36} & 66.95 & 2 & \textbf{4} & -2889.61 & \textbf{-2888.30} \\ \hline 
hotel & 105 & \textbf{78} & 27 & 2.10 & \textbf{1.42} & 6 & \textbf{7} & -4298.20 & \textbf{-4298.11} \\ \hline 
house-dense & 105 & \textbf{92} & 13 & \textbf{2.51} & 3.86 & \textbf{29} & 0 & \textbf{-3778.19} & -3778.52 \\ \hline 
house-sparse & 105 & \textbf{101} & 4 & \textbf{0.12} & 0.16 & \textbf{0} & \textbf{0} & \textbf{-66.78} & \textbf{-66.78} \\ \hline 
motor & 20 & \textbf{13} & 7 & \textbf{1.34} & 1.55 & \textbf{3} & 2 & \textbf{-63.35} & \textbf{-63.35} \\ \hline 
opengm & 4 & \textbf{3} & 1 & \textbf{2.32} & 3.84 & \textbf{2} & 1 & \textbf{-177.17} & \textbf{-177.17} \\ \hline 
pairs & 16 & 3 & \textbf{13} & 1686.86 & \textbf{1487.14} & 3 & \textbf{13} & -68368.80 & \textbf{-68366.13} \\ \hline 
worms & 30 & 0 & \textbf{30} & 4802.06 & \textbf{552.58} & 0 & \textbf{30} & -48504.15 & \textbf{-48497.95} \\ \hline \hline
bur & 8 & 1 & \textbf{7} & 11.99 & \textbf{3.53} & 1 & \textbf{7} & 253431.98 & \textbf{253447.92} \\ \hline 
chr & 14 & 6 & \textbf{8} & 4.95 & \textbf{1.62} & 4 & \textbf{7} & 6097.82 & \textbf{6112.40} \\ \hline 
els & 1 & 0 & \textbf{1} & \textbf{0.02} & \textbf{0.02} & \textbf{0} & \textbf{0} & \textbf{0.00} & \textbf{0.00} \\ \hline 
esc & 18 & \textbf{15} & 3 & \textbf{0.01} & \textbf{0.01} & \textbf{0} & \textbf{0} & \textbf{0.00} & \textbf{0.00} \\ \hline 
had & 5 & 2 & \textbf{3} & 2.27 & \textbf{1.40} & 2 & \textbf{3} & 1556.29 & \textbf{1556.32} \\ \hline 
kra & 3 & \textbf{3} & 0 & \textbf{0.04} & \textbf{0.04} & \textbf{0} & \textbf{0} & \textbf{12133.33} & \textbf{12133.33} \\ \hline 
lipa & 16 & 4 & \textbf{12} & 100.86 & \textbf{59.38} & 4 & \textbf{12} & \textbf{402958.29} & 402886.94 \\ \hline 
nug & 15 & \textbf{12} & 3 & \textbf{0.02} & \textbf{0.02} & \textbf{0} & \textbf{0} & \textbf{244.80} & \textbf{244.80} \\ \hline 
rou & 3 & 0 & \textbf{3} & 4.38 & \textbf{0.44} & 0 & \textbf{3} & 54437.97 & \textbf{54523.46} \\ \hline 
scr & 3 & \textbf{3} & 0 & \textbf{0.00} & 0.01 & \textbf{0} & \textbf{0} & \textbf{46761.33} & \textbf{46761.33} \\ \hline 
sko & 13 & 4 & \textbf{9} & 8.53 & \textbf{3.63} & \textbf{0} & \textbf{0} & \textbf{0.00} & \textbf{0.00} \\ \hline 
ste & 3 & \textbf{2} & 1 & 0.13 & \textbf{0.09} & \textbf{0} & \textbf{0} & \textbf{0.00} & \textbf{0.00} \\ \hline 
tai & 26 & \textbf{13} & \textbf{13} & 13.67 & \textbf{12.23} & 5 & \textbf{10} & 93710.44 & \textbf{93771.96} \\ \hline 
tho & 2 & \textbf{2} & 0 & \textbf{0.08} & 0.09 & \textbf{0} & \textbf{0} & \textbf{0.00} & \textbf{0.00} \\ \hline 
wil & 2 & 0 & \textbf{2} & 6.18 & \textbf{4.18} & \textbf{0} & \textbf{0} & \textbf{0.00} & \textbf{0.00} \\ \hline 
 \end{tabular}
    \caption{The same as Table~\ref{ta:additional_comparison_hungalgs2}, but \emph{worse} bound and time has been computed based on the stopping condition described in Section~\ref{se:compare_attainable_bounds}. See Section~\ref{se:detailed_comparison} for more details.}
    \label{ta:additional_comparison_hungalgs}
\end{table}

\end{landscape}}

\section{Conclusion}

Our theoretical results provide a characterization of the relative interior of the set of optimal solutions of the LP formulation of the LAP. Using this characterization, we were able to provide a linear-time algorithm to calculate such a solution (Algorithm~\ref{al:ri-solution}) and extend these results to the case of ILAP (Section~\ref{se:reduction}).

We employed the aforementioned results in an iterative method that computes a bound on the optimal value of the IQAP. Based on our experiments, our method with relative-interior solution (\hungrialg/) is frequently capable of providing bounds near (or equal to) the optimum of the LP relaxation and is also much faster when compared to a commercial LP solver. The method also typically provided better bounds than the one not following the relative-interior rule (\hungalg/), which confirms the theoretical result from~\cite{werner2020relative}. However, the converse sometimes happened too, i.e., in several instances \hungalg/ was able to achieve a better bound than \hungrialg/. This is expected as each of the methods has many fixed points with different objective values. In practice, the relative performance of these methods may vary depending on the types of instances and on the optimal trade-off between bound quality and runtime.

Concerning future work, our approach might be more useful in pruning the search space during branch-and-bound search where exact solution is sought thanks to the higher quality of the computed bounds. This may apply not only to the (I)QAP, but also to other problems where the LAP naturally occurs as a subproblem, such as the travelling salesperson problem. Finally, although we used the Hungarian method for computing an optimal solution to the LAP, there is large potential to make this computation faster using Sinkhorn~\cite{cuturi2013sinkhorn} or auction~\cite{bertsekas1993reverse} algorithms.

\vspace{1cm} 

\noindent
\textbf{Acknowledgement}

\noindent
Tom\'{a}\v{s} Dlask was supported by the Grant Agency of the Czech Technical University in Prague (grant SGS22/061/OHK3/1T/13), the Czech Science Foundation (grant 19-09967S), the CTU institutional support (future fund), and the OP VVV project CZ.02.1.01/0.0/0.0/16\textunderscore019/0000765. Part of this research was done while Tom\'{a}\v{s} Dlask was visiting the Heidelberg University. Bogdan Savchynskyy was supported by the German Research Foundation (project numbers 498181230 and 539435352).

\bibliographystyle{unsrt}
\bibliography{mybibliography}

\end{document}